\documentclass[a4paper,12pt,leqno]{article}
\usepackage{exscale}
\usepackage{amsmath}
\usepackage{amsfonts}
\usepackage{stmaryrd}
\usepackage{amscd}
\usepackage{amssymb}
\usepackage{theorem}
\usepackage{latexsym}
\usepackage[final]{epsfig}
\usepackage{epsf}
\usepackage{graphicx}

\newfont{\abc}{cmtt10 scaled 1200}

\def\R{\mathbb{R}}
\def\Z{\mathbb{Z}}

\def\ve{\varepsilon}

\def\ra{\rightarrow}
\def\cs{\symbol{35}}
\def\p{\partial}
\def\qed{\hfill $\Box$ \\}

\def\op{\operatorname}
\def\Ric{\op{Ric}}

\def\op{\operatorname}
\def\Ric{\op{Ric}}

\def\supp{\op{supp}}

\def\op{\operatorname}
\def\Ric{\op{Ric}}
\def\scal{\op{scal}}

\def\supp{\op{supp}}
\def\x1n-1{x_1,\dots,x_{n-1}}

\def\dipsi1{\frac{\partial\psi}{\partial x_1}}

\def\d2phin{\frac{\partial^2\psi}{\partial x_n^2}}

\def\1d{\frac{1}{d}}
\def\sx{scal (g) > 0}

\begin{document}

\begin{center}\Large{\bf{Singular Minimal Hypersurfaces and \\Scalar Curvature }}\\
\medskip
{\small{by}}\\
\medskip
\large{\bf{Ulrich Christ \& Joachim Lohkamp}} \\

\end{center}

\vspace{1.7cm}
\noindent Mathematisches Institut, Universit\"at M\"unster, Einsteinstrasse 62, Germany\\
 {\small{\emph{e-mail: uchrist@math.uni-muenster.de, lohkampj@math.uni-muenster.de}}}
\vspace{1.7cm}

\large \textbf {1. Introduction}\\

\bigskip
\normalsize Finding obstructions to positive scalar curvature and getting structural insight is presently based on two competing approaches: one path which is
most travelled  works in the context of spin geometry and gives quite a direct link to topology (cf. [GL1-2] and [G]). The second, much less used but a priori
more general method of attack analyzes minimal hypersurfaces within the manifold under consideration (cf. [SY1-3] and [S]). Although applicable without any
additional (topological, e.g. spin) assumptions and despite its
 natural sensitivity to geometry this approach has a deterrent effect because being based on geometric measure (and regularity) theory it is bound to run into trouble in dimensions above $8$: the appearance of rather hard to understand singularities made the usage of minimal hypersurfaces for
studying scalar
curvature in dimensions $> 8$ basically impracticable.\\

Our aim is to explain how to bypass this problem without losing the information encoded in the singular hypersurface. Conceptionally speaking this is based on
some type of regularizations which are coarse in the sense that the results are not analytically close to the original one but they are fine in the world of
scalar curvature. The central construction in this paper is deleting a carefully chosen neighborhood of the singular set and doubling the resulting manifold in
such a
way that the obtained objects form a sufficiently good substitute for the original hypersurface for use in scalar curvature geometry.\\

In order to state the main result of this paper we consider a closed $n+1$-dimensional Riemannian manifold $(M^{n+1}, g) , n \ge 7$ with \emph{positive} scalar
curvature ($\sx$), a given homology class $\alpha \in H_n (M, \Z)$.\\ Classical geometric measure theory (cf. [D], [F1], [Gi]) provides us with an area
minimizing hypersurface $H^n$ in $(M^{n+1}, g)$ representing $\alpha$ which in general (and of course we assume
this is the case) contains a compact singular set $\Sigma^{n-7}$ of Hausdorff-dimension $\le n-7$, $\emptyset \neq \Sigma^{n-7} \subset H^n$. \\

\textbf{Theorem}{\itshape \quad For any  $\ve > 0$ there is a smoothly bounded neighborhood $V_{\ve} \subset \ve$-neighborhood of $\Sigma^{n-7}$ such that the
 doubling $H^n \setminus  V_{\ve}  \cup_\sim H^n \setminus  V_{\ve} $ admits a smooth metric $g_{\ve}$ with $scal
(g_{\ve}) > 0$. ($\sim$ means gluing along $\p V_{\ve}$)\\

Except for a kind of generalized warped product deformation near $\p V_{\ve}$ the metric $g_{\ve}$ is conformal to the induced metric on $H^n \setminus V_{\ve}
\subset (M^{n+1}, g)$.} \\

Although we heavily use the fact that $\Sigma^{n-7}$ is the set of singularities, the codimension needed to carry out the argument is just $>2$. This matches as
a counterpart of the $Scal > 0$-preserving $codim \ge 3$-surgeries in [GL3] and [SY4] but not along a (tube around some) smooth submanifold within a manifold
with $Scal > 0$:\\ Instead we make up an analytic kind of \emph{stratified surgery} (along submanifolds with positive mean curvature $N^{n-2} \subset H^n$
surrounding strata which in turn will be approximations for $\Sigma^{n-7}$) in a space whose \emph{first eigenvalue} for
(a scaling invariant refinement of) the \emph{conformal Laplacian is positive}.\\

 The techniques described in this paper can be used and extended to handle obstruction theory for $scal > 0$ in \emph{arbitrary dimensions}
 non-existence of $Scal > 0$-metrics on enlargeable manifolds (e.g. $T^n \cs N^n$)[L1], more general lower scalar curvature bounds.\\
For non-compact manifolds one gets results for sufficiently tame ends, i.e. product like or asymptotically flat resp. hyperbolic ends. In particular, we can use
this result to derive short geometric proofs of the general positive mass conjectures in every dimension [L2] which extend to more advanced versions (e.g. with
certain non-asymptotically flat complete ends)

Moreover, we also note that, for instance, $Spin^{\mathbb{C}}$-problems  can be shifted to $Spin$-geometry (after taking a suitable $S^1$-bundle) while loosing
some scalar curvature information (cf. [LM]) which can partially be recovered taking minimal hypersurfaces. Thus the freedom to extend the dimensional range
also provides additional tools for lower dimensions.\\

Since the proof is a bit involved (although we think that the ideas are quite natural) we describe it explaining roughly the geometric effects (A), the chronology of the arguments (B) and a conceptional viewpoint (C).\\

For (A) we can think of three main steps: \\
\textbf{1}. We deform the singular hypersurface $H^n$ outside its singular set $\Sigma^{n-7}$ into some pointwise $\sx$-geometry. (Remark: \emph{n-7} is just an
upper bound but the dimension may globally but also locally be smaller -- writing  $\Sigma^{n-7}$
refers just to this upper dimensional bound but \emph{not} to its actual dimension.)\\
After spending some time on deriving estimates for the new metric in particular relative to the original one we will be able to carry
out deformations along  $\Sigma^{n-7}$. \\
This \textbf{2}nd step can be thought of as a stratified version of  $codim \ge 3$ surgery for positive scalar curvature along an (actually augmented) singular
set. In the classical regular case one gets a totally geodesic boundary keeping the scalar curvature $>0$. The counterpart we obtain is an implicit
barrier for $(n-2)$-dimensional minimal hypersurfaces $\subset H^n$ homologically equivalent to the boundary of a neighborhood of $\Sigma^{n-7}$.\\
\textbf{ 3}. From this we get a \emph{smooth} $(n-2)$-dimensional hypersurface $N^{n-2}$ with positive mean curvature homologically equivalent to that boundary and
arbitrarily close to $\Sigma^{n-7}$. Now a non-conformal deformation transforms a small one sided tube of $N^{n-2}$ into a totally geodesic border (and additionally gives some extra $\sx$). Gluing this with a mirrored copy completes the argument.\\

Now, for (B), let us give some more details, basics and notations used later on: The deformations in step 1 and 2 are conformal. They transform the  metric
$g_H$ induced on H from $(M^{n+1}, g)$ to $u^{4/n-2} \cdot g_H$ (for some smooth $u > 0$ defined on $H \setminus \Sigma$) and the transformation law (TL) for
the new scalar curvature is
\[   - 4 (n-1)/(n-2) \cdot \Delta u + scal_{g_H} \cdot u =
 scal_{u^{4/n-2} \cdot g_H} \cdot u^{n+2/n-2} \]

Specifically, the start point for the first step is the observation by Schoen and Yau that in the case where $H$ is a smooth closed hypersurface the fast that
$H$ is an area minimizer
 implies that the first eigenvalue of the left hand side operator $L$ in (TL) (the conformal Laplacian) $L u = - \frac{4(n-1)}{n-2} \cdot \Delta u + Scal_H \cdot u$
has a \emph{positive} first eigenvalue $\lambda_1$:\\ Since $H$ is area minimizing the 2nd variation of its area is $\ge 0$, that is $H$ is a \emph{stable}
minimal hypersurface. Formally, $A''(f) := Area''(f \cdot \nu) \ge 0$ where $\nu$ is a unit normal vector field (we may assume $M$ and $H$ are orientable) and
$f \in C^\infty \mbox{ or } H^{1,2} (M, \R)$, thus $f \cdot \nu$ is an infinitesimal variation of $H$ and a direct computation gives the expression:
\[ (A2) \;\;\; A'' (f) = \int_H | \nabla_H f |^2 - f^2 (| A |^2 + Ric_M (\nu, v)) d A \ge 0 \; \Longleftrightarrow \]
{ \small\[ \int_H | \nabla f |^2 + \frac{n-2}{4 (n-1)} scal_H f^2 d A \ge \int_H \frac{n}{2 (n-1)} | \nabla f |^2 + \frac{n- 2}{2 (n-1)} f^2 \left( | A |^2 +
scal_M \right) d A\] }

 where $| A
|^2 = \Sigma^{n-1}_{i=1} a_i^2, a_i = i-th$ principal curvature of $A$ and $\nu$ is the
unit normal vector to $H$. For the equivalence one just uses the identity\\
$ | A |^2 + Ric_M (\nu, \nu) = \frac{1}{2} \left( scal_M - scal_H + (tr A)^2 + | A |^2 \right)$.\\

 Thus if $(M, g)$ has $\sx$ we
get $ \int_H | \nabla f |^2 + \frac{n-2}{4 (n-1)} scal_H f^2 d A \ge \int_H c \cdot f^2$ for $c = \frac{n- 2}{2 (n-1)} \cdot \inf_M scal(g_M) > 0$ and every
smooth function $f$ on $H$ and therefore \[ \lambda_1 = \inf \{ \int_H | \nabla f |^2 + \frac{n-2}{4 (n-1)} scal_H f^2 d A \;|\; f \in H^{1,2}(M), |f|_{L^2} = 1 \}
\ge  c > 0. \]

Finally, using the fact that the first eigenfunction $u_1$ can be assumed to be positive we observe that $scal_{u_1^{4/n-2} \cdot g_H} > 0$, namely:
\[ \lambda_1 \cdot u_1 = - 4 (n-1)/(n-2) \cdot \Delta u_1 + scal_{g_H} \cdot u_1 =
 scal_{u_1^{4/n-2} \cdot g_H} \cdot u_1^{n+2/n-2} \]

Now we switch to the case where $H$ has a \emph{non-trivial singular set} $\Sigma$. It is a classical result (cf. [D],[Gi] and [Si]) that this \emph{compact}
set has at least \emph{codimension 7} within $H^n$. In dimension $8$ there are only isolated singular points and in this particular situation one has a theory
(cf.[HS], [Sm]) of $C^k$-perturbing the metric to resolve the singularities.
 But in higher dimensions the structures is more or less unknown (even rectifiability is unclear) and corresponding perturbation results are out of reach.\\
The singular set could be a fractal set and will usually have components of varying Hausdorff-dimension $\le n-7$ in $M^{n+1}$.
 However, as a by-product, we will construct a series of \emph{stratified spaces}
approaching $\Sigma$ in an analytically effective way. \\

 In order to handle this case we first note that the area minimizing property of $H$
enters significantly in our discussion of the scalar curvature geometry near the singular set, whereas the smooth case relies only on the \emph{stability of}
$H$: More specifically we use the fact that $codim (\Sigma)
> 2$ for the original hypersurface (in turn this - when supposed as an extra data - already implies $codim (\Sigma) \ge 7$ for stable hypersurfaces acc.[SSi]) but also
for its tangent cones \emph{and} we use that this property singles out a \emph{compact} (not just precompact) set in the set of minimal hypersurfaces respectively  minimal cones. \\

Nevertheless, we first use that $A''(f) \ge 0$ for all functions $f$ with $supp(f) \subset H \setminus \Sigma$, but in contrast to the smooth case we use the
particular entities of (A2) in a substantial way:\\ The term $\int_H  \frac{n- 2}{2 (n-1)} f^2 | A |^2  d A$ provides us with the needed stronger estimates. In
particular, the first eigenvalue of an asymptotically "scaling invariant" form of the conformal Laplacian is still positive. Specifically we work with
\[ - \triangle u_0 + \frac{n-2}{4 (n-1)}  \scal_H u_0 = \lambda_0 \cdot  | A | ^2  \cdot u_0\]

This is motivated by the idea that with this kind of equations the positivity of the eigenvalue survives the passage to tangent cones approximating $H$ in
$\Sigma$ where we get a better (more precisely:
inductive) control over the behavior (note that (under scaling) the eigenvalue of the conformal Laplacian would converge to zero and thus we would lose the
key information in this transition).\\

This allows us to become a more specific about the deformation in the Theorem not too close to $\Sigma$:\\

\textbf{Addendum }{\itshape \quad On $H^n \setminus  V_{\ve} $ the metric $g_{\ve}$ is $u_0^{4/n-2} \cdot g_H$ up to some conformal redistruction of $scal >
0$.}\\

Next we construct (rather than get) an eigenfunction for $\rho^2 L$: due to the ill-posedness of this problem the function will not give any useful insight when
approaching $\Sigma$. Here we work out a technique to characterize a particular (super-)solution modifying a Perron-type construction of solutions of certain
types of elliptic equations. The main feature of these particular solutions which are in a sense \emph{minimal} is that there defining properties are
inherited by the tangent cones and due to that scaling invariance we find the same type of equation and solutions on the cones.\\
In turn the limiting behavior of this supersolution $w$ near $\Sigma$ can be captured using transition to tangent cones where these solutions admit
a separation of variables and hence become easier to understand, in particular they are amenable to an inductive scheme modelled over the classical dimension reduction for
area minimizers cf. [F2]. However the justification of this method is quite involved and a subject on its own that appears in [L3]. \\

To proceed with the argument one noticed an interesting effect induced from the Perron-type construction: close to points $p \in \Sigma$ the deformed metric
still looks like a cone: $(C,\tilde g)$ is isometric to any of copy scaled around $0$ and can be reparametrized as $c(\omega)^{4/n-2}  \cdot g_{\R} + r^2 \cdot
g_{\p B_1(0) \cap C}$. A tricky point will be that this closeness depends
\emph{discontinuously}  on the base point $p \in \Sigma$.\\

Next we analyze the resulting metric $w^{4/n-2} \cdot g_H$ and find that we can define some kind of Greens functions on small balls and inductively on strata of
an augmentation of $\Sigma$ whose shape and boundary behavior can also be controlled.\\
Here the fact that the "Hausdorff codimension" of $\Sigma$ is $>2$ becomes critical. A non-constructive covering and localization strategy from [L5] allows us
to find a finite collection of such functions defined on a collection of balls and strata covering $\Sigma$ according to the definition of the Hausdorff measure
is added and used to subsequently deform $(H, w^{4/n-2} \cdot g_H)$ along $\Sigma$ keeping $\sx$ but producing a barrier such that area minimizing hypersurfaces
in $N_0^{n-1} \subset H^n$ homologically equivalent to a boundary of a
neighborhood of $\Sigma$ are deflected keeping them away from $\Sigma$. \\

If we now try to find such an area minimizing hypersurface $N_0^{n-1} \subset H^n$  we will definitely need and want to restrict the region close to $\Sigma$
where the support
of $H$ is supposed to be. Actually, the absolute area minimizer in this homology class would be a point.\\

In order to accomplish this (and for other more technical purposes as well) we introduce \emph{parametric minimal hypersurfaces with obstacles} designed to
smooth this minimizer: we place a separating collection of balls outside a neighborhood of $\Sigma$. While an area minimizer $N^{n-1}$ in this region just no
longer reaches $\Sigma$ it touches some of these balls which form a rigid barrier. The intersection sets are harmless but the rest of $N^{n-1}$ is just a free
minimal hypersurface which may have the usual singularities. However after placing the balls suitably (which is a rather implicit construction) we conclude that
the resulting hypersurface is $C^1$. Here we use the boundary regularity result of Allard to conclude that sufficient control of the boundary of the
intersection of $N^{n-1}$ with the barriers provides us with estimates for the distance
to the first potentially singular point. The hypersurface $N^{n-1}$ can then be smoothed with resulting mean curvature $ \ge 0$. \\
Hence a non-conformal deformation which compresses the geodesics leaving the hypersurface perpendicularly directed towards $\Sigma$ additionally increases
the scalar curvature and gives a totally geodesic boundary  pr\^ et-\`a-porter for gluing with a mirror copy.\\

Finally, for (C), a conceptional remark. The techniques are designed to bundle the data of the area minimizer $H$, certain uniquely characterized functions and
submanifolds on/in $H$ etc. to canonical objects which obey \emph{compactness theorems}. Then one can consider limit objects (usually area minimizing cones) and
reduce the estimates and certain constructions by inductive cone reduction arguments (which are extensions of the classical way to roughly analyze $\Sigma$ by
successive blow ups (cf. [F2], [Gi],Ch.11)) until one reaches dimension $8$.  Now one reverses the direction and assembles a geometry on $H$ inductively by
those covering argument mentioned above. \\

\vspace{1.7cm}
\large \textbf {2. Area Minimizing Cones and Reduction Techniques}\\
\normalsize

The only a priori information concerning the singular set $\Sigma \subset H^n$ we use is the compactness and the Hausdorff-dimension which is $\le n-7$.
 But we have a structural aid provided by tangent cones (cf. [Gi], [Si]). These are (locally area minimizing) minimal cones in
$\R^n$ forming a generalization of the tangent plane at regular points: after some scaling one may consider $H$ locally (say around $p \in \Sigma$)
being embedded in $\R^n$ and after further scalings by an increasing sequence of factors ($\tau_m \ra +\infty$) there is a minimal cone $C_p$ which approximates
$\tau_m \cdot H$
on any given compact set in $\R^n$ in a certain way described below. The point is that (in sense we will have to discuss) this is a linearization in one direction.\\
The usage of tangent cones in the literature is fairly limited since each singular point in $\Sigma \subset H$ will usually have infinitely  many tangent cones,
the set of tangent cones varies discontinuously alon g $\Sigma$ and the approximation of $H$ by these cones is {\bf not} uniform in $\Sigma$. \\
Nevertheless, in part because we will be able to avoid to come \emph{too} close to $\Sigma$, we can set up a scheme to derive many properties of $H$ near
$\Sigma$ from corresponding information on cones and for certain properties this allows us to gain \emph{uniform} control using the various cones as a link. \\
And, more importantly, we can carry out certain local operations of cones (serving as models) and transplant to $H$.\\

We start on an abstract level with a composition of several classical facts due to De Giorgi, Allard and
others (cf. [D], [A1] , [Gi] and [Si]). \\

{\bf Proposition (2.1)}{\itshape \quad  Let $H^n \subset M^{n+1}$ be an area minimizing hypersurface and $\tau_m \to +\infty$ a sequence of positive real
numbers.\\ Then, for every $p \in \Sigma$ we find a subsequence $\tau_{m_k}$ and a cone $C_p \subset \R^n$ such that for any given open $U \subset \R^n$ with
compact closure the \textbf{flat norm} $d_U$ (cf. [S],Ch.31) which (roughly speaking) measures the volume between two sets in U  converges to zero:
\[ d_U (\tau_{m_k} \cdot H, C_p) \to 0 \]
and this convergence implies \textbf{compact $C^l$-\emph convergence}, for any $l \ge 0$, if $\overline U$
contains only smooth points of $C_p$.} \\

{\bf Remark} \quad $\tau_{m_k} \cdot H \subset \tau_{m_k} \cdot M$ can locally (near $p$) be considered a subset of $\R^n$ (for $k \to +\infty$ the deviation
vanishes). The $C^l$-convergence statement can be formulated more precisely as follows: assume V is an open subset of $C_p$ whose compact closure contains only
regular points, the focal distance $\iota(V)$ is $> 0$ and we consider the $exp_\nu$-image $U_\epsilon$ of normal vectors of length $\le \epsilon \le \iota/2$
in the normal bundle $\nu|_V$ of $V \subset C_p$. Then for large $k$ the set $U_\epsilon \cap \tau_{m_k} \cdot H$ is a $C^l$-graph (= $C^l$-section of the
normal bundle) over
 $V$ and converges compactly to $V$ (= zero section) in $C^l$-topology. \\

The cone reduction argument we are looking for cannot be based on particular properties of a special cone but becomes valid only if we can manifest such
properties for the class of all singular cones
simultaneously. One of the ingredients will therefore be the following two results \\

{\bf Lemma (2.2)}{\itshape \quad  The set of embedded area
minimizing $n$-cones $C_{n}$ (around $0$) in
$\R^n$ is compact in the flat norm topology. }\\

{\bf Proof} \quad This can be derived from 37.2 in [Si]. More
directly, one may use the regularity (mod \emph{codim 7}
singularities) to inductively derive an explicit bound for the
$n$-dimensional volume $\cap B_1(0)$ and (hence the
$(n-1)$-dimensional volume $\cap \p B_1(0)$) by some comparison
with the unit sphere $\p B_1(0)$. The compactness theorem for
integral currents and the fact that minimality and the cone shape
survive under flat norm convergence (e.g.
using Allard regularity for convenience) give the result. \qed \\

In particular, the set  ${\cal T}_H$ of tangent cones of $H$ (with center set to 0) has the \emph{compact} closure $\overline {\cal T}_H \subset C_{n}$. $\p
\overline {\cal T}_H$ will usually contain cones which do {\bf not} appear as tangent cones of $H$. Actually considering such extensions deliberatively will be
an essential tool for many arguments: we
state a simple but crucial compactness result in this direction: \\

{\bf Corollary (2.3)} {\itshape \quad There is a constant $d_n >0$ such that \\

$d_{B_1(0) \setminus B_{1/2}(0)} (C, \R^{n}) < d_n$ \; if and only if \; $C$ is non-singular. \\

Therefore the set of singular
cones  $SC_{n} \subset C_{n}$ is closed (and hence compact) and hence $\overline {\cal T}_H \subset SC_{n} $.} \\

{\bf Proof} \quad Let $C_i$ be a sequence of such cones with $d_{B_1(0) \setminus B_{1/2}(0)} (C, \R^{n}) \to 0$; then by the cone property
$d_{B_1(0)}(C,\R^{n}) \to 0$ and Allard
regularity implies that for large $i$ every $C_i$ is non-singular. \qed \\

Later on we will derive several other universal properties of singular cones using Allard and the .\\

Next we will sharpen the usual picture of cone approximation: For decreasing radius $\eta \ra 0$ $(\eta^{-2} \cdot H) \cap B_2(p) \setminus B_1 (p)$ is not just
sometimes approximated by a cone but a slightly closer look already unveils an instructive view: choose a
\emph{finite} covering $\{ B_\delta (c_i) \}$ of the compact set of singular cones $\cal{CS}$ by \emph{flat norm} balls of radius $\delta$.\\
1. For any $\delta
> 0$ we find that starting from some $\eta_\delta
> 0$ such that $(\eta^{-2} \cdot H) \cap B_2(p) \setminus B_1 (p)$ is $\delta > 0$ - close in flat norm to some (non uniquely determined) tangent cone
$C^\eta_p$.\\
 2. Considering this assignment as a discrete valued map $\eta \mapsto \{ B_\delta (c_i) \}$ we observe a large scale \emph{fading} or \emph{freezing} property:
 after scaling $\eta$ to 1 the frequency of oscillation within
the balls of this finite covering will decay uniformly to zero for  $\eta \ra 0$ and (also after scaling) the size of the well-approximated part
of \emph{any} of these cones increases (i.e. considering a sequence of approximating regions (identified via scaling) we get a compact exhaustion of any tangent cone).\\

This is just an interpretation of the following \\

{\bf Lemma (2.4)} \quad \emph{  For any $\delta > 0$ and any pair $R \gg 1 \gg \varrho >0$ we can find a small $\eta_{\delta , R , \varrho} > 0$ such that for
\textbf{every} $\eta \in (0, \eta_{\delta , R , \varrho})$:
\[\eta^{-2} \cdot H \cap B_R(p) \setminus B_\varrho (p) \mbox{ is }\delta > 0 - \mbox{close in flat norm to a tangent cone } C^\eta_p.\]}

(From Allard regularity one gets corresponding statements in terms of $C^k$-topology on smooth parts of $C_p$. Note also that $\eta_{\delta , R , \varrho}$
depends on $p$ in an
discontinuous way. )\\

The \textbf{proof} is standard: if there were a sequence of $\eta_i \ra 0$ and a $\delta_0 > 0$ such that $\eta^{-2} \cdot H \cap B_R(p) \setminus B_\varrho
(p)$ is not  $\delta_0$-close to any tangent cone, there is still a subsequence that eventually being arbitrarily close to a tangent cone leading to an
immediate contradiction.\qed \\

Notice that Allard regularity provides us with the refined version for $C^k$-topology. We set for any tangent cone $V_{\xi}(\sigma) = \mbox{\emph{cone over}
}U_{\xi}(\sigma \cap \p B_1(0)) \subset  C^\eta_p$ such a smooth scaled version $H$: $\eta^{-2} \cdot H \cap B_R(p) \setminus B_\varrho (p) \setminus
V_{\xi}(\sigma) $  can be written as a graph of a function $g_\eta$ over
a corresponding part of $C^\eta_p$ (in the sense explained above). Now (2.4) and (2.1) imply:\\

 {\bf Corollary (2.5)} \quad \emph{  For any $\delta > 0$ and any triple $R \gg 1 \gg \varrho \gg \xi >0$ we can find a small $\eta_{\delta , R , \varrho, \xi}
> 0$ such that for \textbf{every} $\eta \in (0, \eta_{\delta , R , \varrho, \xi})$:
\[ |g_\eta|_{C^k} < \delta \mbox{ on } (B_{R}(0) \setminus B_{\varrho}(0)) \setminus V_{\xi}(\sigma) \subset  C^\eta_p\] } \qed

We now describe the basic procedure we use to mediate between $H$ and the realm of singular cones and
how to proceed from there. \\

In order to prove a local result on $H$ which is known to be true for cones we permanently argue by contradiction. Assume there is a sequence of points $x_n \in
H \setminus \Sigma$, ${\rm dist}_M(x_n, \Sigma) = \ve_n \to 0$ (we have to consider the intrinsic distance later on) and around $x_n$ a certain expected
geometric (or more general analytic) property fails to hold on $B_{\alpha \cdot \ve_n}(x_n)$, $\alpha \ll 1$. In addition, the property in question should satisfy a compactness property:
elliptic compactness (and Arzela-Ascoli) when we consider eigenfunctions, Gromov compactness (plus Allard
regularity) when we consider the second fundamental form as a curvature quantity. \\
We will then argue as follows: There is a $p \in \Sigma$, being limit of a subsequence of $x_n$ and $\rho_n = d(x_n,p) \ge \ve_n$ will also converge to zero.
However, there are two cases
\begin{enumerate}
\item $\ve_n / \rho_n > const. > 0$, in this case the $x_n$ run into a well-approximated zone of a tangent cone in $p$ \item $\ve_n / \rho_n \ra 0$, here we
still get a cone approximation, but the cone may not appear as a tangent cone in any point
\end{enumerate}
In case (i), after scaling $H$ and $M$ by $\ve_n^{-2}$, since $d (x_n,p)$ is now normalized to 1 (up to bounded multiple), there is still a subsequence of
$x_{n_k}$ converging (in this scaled picture) to a point $q \in \p B_1(0) \cap
C_p$ where $C_p$ is a tangent cone at $p$. \\

Now in case (ii), we can argue as follows: Take a point $p_n \in \Sigma$ with $d(x_n,p_n) = dist (x_n, \Sigma) =  \ve_n$  and scale each intersection $H \cap
B_{ \rho_n}(p_n)$ by $\rho_n^{-2}$. This can be considered a sequence of area minimizing surfaces $T_n$ in $B_1(0) \subset \R^n$ and we may assume it converges
in flat norm to an area minimizer $T_\infty$ in $B_1(0) \subset \R^n$. A subsequence of $(\ve_n / \rho_n)^{-2}$-scaled copies of $T_\infty$ converges in flat
norm to
a minimal cone $C_\infty$ (which may not be a tangent cone of $H$). \\
Thus, by a diagonal sequence argument we may assume that   $H \cap B_{ \rho_n}(p_n)$ scaled by $\ve_n^{-2}$ converges in flat norm to $C_\infty$ and that
 $x_{n}$ converges (in this scaled picture) to a point $q \in \p B_1(0) \cap C_\infty$. \\

 In both cases the limiting cone is smooth outside a \emph{codim 7} singular set $\sigma$, hence the flat norm convergence gives rise to compact
$C^l$-convergence
outside $\sigma$. \\
For  convenience we will use $C_\ast$ as a common notation for $C_p$ resp. $C_\infty$ when both cases can show up. The second case will also be called an
\emph{abstract cone reduction}.\\

The cone reduction strategy proceeds as follows: In certain cases an \emph{a posteriori} argument shows that $q$ is a \emph{regular point} in  $C_\ast$. In some
other cases we use that after scaling around $q$, $C_\ast$ can be approximated by a tangent cone which is a product $\R \times \hat{C}^{n}$, where $\hat{C}^{n}
\subset \R^{n}$ is again a minimal cone and argue inductively. Now we may use the compactness result for the geometric/analytic estimate or property under
consideration and the fact that $\ve_n^{-2} \cdot H$ converges to $C_\ast$ to conclude the estimate/property continues to fail on $B_\alpha(q)
\subset C_\ast$. Therefore we are done if we know that, in fact, the corresponding property does hold on $B_\alpha(q)$. \\

Direct arguments (and hence sharper estimates) often fail since this would usually require uniform approximation by
tangent cones. \\

The cone structure actually provides us with two tools: the cone direction which blows up to give a local product structure with a minimal hypersurface
$G^{n-1}$ is used as a construction aid on its own and, secondly, the properties of $G^{n-1}$ can be used as
induction hypothesis. \\

{\bf Remark} \quad At this point it is important to mention that $G^{n-1} = \p B_1(0) \cap C$ is minimal but neither area minimizing nor stable (since $\p
B_1(0)$ has $\Ric
>0$). Nevertheless we can carry over those results valid for area minimizers which allow us to make the induction work: The crucial property of $G^{n-1}$ in
this setting is that the \emph{cone over} $G^{n-1}$ is area minimizing and therefore all its tangent cones are. Outside 0 the tangent cones have a Riemannian
product structure isometric to $\R \times \tilde C_q^{n-1}$ where $\tilde C_q^{n-1}$ is again area minimizing. However, theses cones $\tilde C_q^{n-1}$ are
precisely the tangent cones of $G^{n-1}$. This, together with the local product structure of $C$ as a cone over $G^{n-1}$ will allow us to handle $G^{n-1}$ is
our scheme just like an actual area minimizer, e.g. the singular set of $G^{n-1}$ has the same properties (e.g. ${\rm codim} \ge 7, compactness)$ as for area
minimizers.
In addition, the argument for the two cases of distinction for $C_\ast$ survive. This would not be the case for general minimal surfaces. \\

Also, there will be no accumulating problem during the induction process, since the next step is to stick with the tangent cones $\tilde C_q^{n-1}$ of
$G^{n-1}$, consider $\p B_1^{n-1}(0) \subset \R^{n-1}$ and $G^{n-2}$ until we obtain isolated point
singularities. \\
A good way of thinking of this part of the strategy is as an enhancement of the classical cone reduction in determining the codimension of $\Sigma$ with
additional data on the hypersurfaces which induce corresponding data on the lower dimensional objects.\\

As a sample of this rather abstract scheme we consider the \emph{intrinsic} distance function on $H$. The distance between points $x \in H \setminus \Sigma$ and
(points in) the compact set $\Sigma \subset M$ measured within the ambient manifold (extrinsic distance) ($d_M(x,p)$ resp.) ${dist}_M(x, \Sigma)$ is not
suitable for our purposes: We use the intrinsic metric on $H$ to study e.g. eigenfunctions of the conformal Laplacian.  Also we will conformally deform the
induced metric on $H$ and thereafter want to understand the behavior of the new geometry near
$\Sigma$, but at that stage the embedding has lost its meaning. Thus we have to work with the intrinsic distance function ${dist}_H(x, \Sigma)$ on $H$. \\
Since $H$ may develop additional bumps and even new topology when approaching $\Sigma$ (reflected by thin regions with large $|A|$) one realizes that it is by
not at all clear that $d_H(p,x) < +\infty$ for any two points $x \in H^n \setminus \Sigma$, $p \in \Sigma$ in the same connected component of $H$.\\ However,
using the fact that $H$ is an area minimizer we can actually prove below that close to $p$ there is a network of pieces of rays which links to
$p$ in finite time. We will base this argument on a cone reduction. \\

We first note a valuable relation "$\Sigma \prec \sigma$" (which is (2.6) below) between the singular sets $\Sigma \subset H$ and $\sigma \subset C_p$ :
asymptotically the singular set of the tangent cones is  "larger" than the germ of the singular set around $p \in H$. For instance, $\Sigma$ may contain
scattered points or there might be smooth but highly curved regions near $\Sigma$ which may cause the appearance of asymptotically tangent rays in $\sigma$ . On
the other hand, the complexity of
the $\sigma$ is reduced by one dimension (since $\sigma$ is also a cone). \\

Moreover the critical deformations close to $\Sigma$ which will be handled later on will be prepared on cones (instead of $H$) and will be transplanted to balls
in $H$ keeping entirely their effect because the resulting "horizon" will hide away still smaller neighborhoods of $\sigma$ and, thus by "$\Sigma \prec \sigma$"
we are
able to imitate a stratification of $\Sigma$ by local enhancements of $\Sigma$ by $\sigma$:\\

{\bf Lemma (2.6)}: {\itshape \quad Let $U \subset \R^n$ be an arbitrarily small neighborhood of $\sigma \subset C_p \subset \R^n$. Then, for each sequence $t_m
\to \infty$, there is a subsequence $t_{m_k}$ such that
\[ d_{(B_R(0) \setminus B_r(0)) \setminus U} (C_p, t_{m_k} \cdot H) \to 0 \]
for any fixed $R>r>0$. Therefore, for sufficiently large $k$,
\[ t_{m_k} \cdot \Sigma \cap (B_R(0) \setminus B_r(0)) \subset U. \]
(We may therefore think of $\Sigma$ locally being already enveloped by a tiny neighborhood of $\sigma$.)} \\

The proof is just a standard application of Allard regularity; for simplicity we will denote the subsequence by $t_m$ again: $t_m \cdot H$ converges in flat
norm to the smooth $C_p$ in $(B_R(0) \setminus B_r(0)) \setminus U)$, thus $t_m \cdot H$ is also smooth for sufficiently large $m$. In other words, no point of
$\Sigma$ can be contained in $(B_R(0) \setminus B_r(0)) \setminus U)$.\\

{\bf Corollary (2.7) }{\itshape \quad
\begin{enumerate}
\item Let $p \in \Sigma$. For sufficiently small $\rho >0$ assume that $B_\rho(p)
    \cap H$ is connected. Then $(B_\rho(p) \cap H) \setminus \Sigma$ is also connected and
     $d_{H \setminus \Sigma} (p,x) \le c \cdot \rho$ for $x \in B_\rho(p)$.
    \item If $H$ is connected, then also $H \setminus \Sigma$ is connected and its
    intrinsic diameter is finite.
\end{enumerate}}\

(The way we prove this is chosen to be able to extend these arguments directly to situations where we deformed $H$ and recover the new distances considering the
induced
geometries on  tangent cones.)\\

{\bf Proof} \quad If $\Sigma$ is a finite set (and all tangent cones are regular with singularities only in 0) then the tangent cones are connected by the
maximum principle and, due to codimension $\ge 2$, removing $\{0\}$ keeps the complement connected. For $\rho >0$ small enough we can assume that (after scaling
with $\rho^{-2}$) the set $(B_2(0) \setminus B_{1/2}(0)) \cap C$ is $C^k$-close to $B_2(0) \setminus B_{1/2}(0) \cap H$ for a suitable tangent cone $C$ (the
choice depending on $\rho$). (Otherwise there would be a sequence $\rho_m \to 0$ for which there is no $C^k$-close tangent cone, which would imply that there is
no flat norm convergence to any tangent cone for this sequence.) Thus in this case $(B_\rho(p) \cap H) \setminus \Sigma$ can be written as a union of connected
sets
(rescaled versions of $B_2(0) \setminus B_{1/2}(0) \cap H$) and is also connected. \\
In this case the claims concerning intrinsic distances are obvious. \\

Now proceed with the case where the tangent cones also contain singularities other than 0. We claim that for a given $p \in \Sigma$ and $\rho >0$ small enough,
$(B_{2\rho}(p) \setminus B_{\rho/4}(p)) \cap H$ contains an open connected subset $V_\rho$ with the following properties:
\[ {\rm diam}_{V_\rho} V_\rho \le k_n \cdot \rho, \quad \frac{{\rm vol} (V_\rho \cap \p B_{2\rho})}{{\rm vol} (H \cap \p B_{2\rho})} > \frac{3}{4}, \quad \frac{{\rm vol} (V_\rho \cap \p B_{\rho/4})}{{\rm vol} (H \cap \p B_{\rho/4})} > \frac{3}{4} \]
(where ${\rm diam}_{V_\rho} V_\rho$ is the diameter measured within $V_\rho$.)
Assume there is a sequence $\rho_n \to 0$ such that
$B_{2\rho_n}(p) \setminus B_{\rho_n/4}(p)$ does not contain such a
subset $V_{\rho_n}$. We may assume that for this sequence there is
a fixed tangent cone $C_p$ (as the limit object). Since
$B_{\rho_n}(p) \cap H$ is assumed to be connected, $\p B_1(0) \cap
C_p$ is connected (a maximum principle argument gives this anyway)
and by induction (cf. the Remark concerning minimizers above): $\p
B_1(0) \cap C_p \setminus \sigma$ is connected, where $\sigma$ is
the singular set of $C_p$, the intrinsic diameter of $\p B_1(0)
\cap C_p$ is finite and there is a connected $W \subset \p B_1(0)
\cap C_p \setminus \sigma$
with $\frac{{\rm vol} (W)}{{\rm vol} (C_p \cap \p B_1(0))} > \frac{4}{5}$, ${\rm diam}_W W < +\infty$. \\

The compactness result for tangent cones allows us to adapt the choices to get a uniform upper bound $b$ for the respective ${\rm diam}_W W$ for all tangent
cones used in the construction: otherwise there is sequence of tangent cones $C(k)$ converging to some $C$ (in $C^l$ on smooth parts) where the infimum of ${\rm
diam}_W W$ diverges. But we can find a connected $W_C \subset \p B_1(0) \cap C \setminus \sigma$ with $\frac{{\rm vol} (W_C)}{{\rm vol} (C \cap \p B_1(0))} >
\frac{9}{10}$, ${\rm diam}_{W_C} W_C < +\infty$ and for large $k$ this provides us (via $C^l$-identification) with some connected subset $W$ in $\p B_1(0) \cap
C(k) \setminus \sigma$,
$\frac{{\rm vol} (W)}{{\rm vol} (C(k) \cap \p B_1(0))} > \frac{4}{5}$, ${\rm diam}_W W < {\rm diam}_{W_C} W_C + 1$.\\

 Thus defining $\tilde V_{\rho_n} =
B_{2\rho_n}(p) \setminus B_{\rho_n/4}(p) \cap \mbox{ subcone of } W \subset C_p$ and using the $C^l$-approximation
of the scaled $H$ this induces a corresponding set in $H$ giving a contradiction.  \\

Now choosing a $V_\rho$ for each small $\rho >0$ consider $\cup_{k=0}^\infty V_{\rho/2^k}$ (note that now there are several tangent cones involved). Because of
the volume fraction $> \frac{3}{4}$ belonging to $V_{\rho/2^k}$ in each  boundary, we have an open non-compact intersection $(\p B_{\rho/2}(p) \cap V_\rho) \cap
(\p B_{\rho/2}(p) \cap V_{\rho/2})$. Starting at $x$ we now choose a path which follows the (approximate) ray direction in $V_\rho \cap (B_\rho(p) \setminus
B_{\rho/2}(p))$. Then, on $\p B_\rho(p)$ one uses ${\rm diam}_W W < b$ to run to a point which, when following the ray direction, leads (within $V_{\rho/2}$) to
an intersection point with $V_{\rho/4}$ etc. Thus we get a sequence of points $x_m \in H \setminus \Sigma$ with $x_m \to p$, $d_{H \setminus \Sigma}(x_m,
x_{x+1}) \le c \cdot 2^{-m} \cdot \rho$ and conclude $d_{H \setminus \Sigma}(p,x) \le c \cdot \rho$ for $x \in B_\rho(p)$ where $c=1+b$. The other claims are
direct consequences of this construction. \qed

{\bf Corollary (2.8) }{\itshape \quad There is are universal bounds $0 < A_1(n) < A_2(n) < \infty$ and $0 < D_1(n) < D_2(n) < \infty$ for the area $A$ and
diameter $diam$ of $\p B_1(0) \cap C$ for any $C \in SC_n$: \[A_1(n) < A < A_2(n) \mbox { and } D_1(n) <  diam <  D_2(n) \] }\
 {\bf Proof} \quad This is a
consequence of the compactness of $SC_n$ and we only indicate the argument for the least obvious claim $diam_{\p B_1(0) \cap C} < D_2(n)$. If $C_j$ is a
sequence with $diam_{\p B_1(0) \cap C_j} \ra \infty$ we may assume it converges in flat norm and $C^k$-compactly to some limit cone $C_\infty$ and (2.7) says
that $diam_{\p B_1(0) \cap C_\infty} =: D < \infty$. The compact $C^k$-convergence implies that there is a sequence $\ve_j \ra 0$ such that $diam_{\p B_1(0)
\cap C_j \cap U_{\ve_j}} \ra \infty$ where $U_{\ve_j}$ is the extrinsic $\ve_j$-neighborhood of $\sigma_j$. Rescaled by $\ve^{-2}_j$ we reach for large $j$ (cf.
remark after (2.5)) a minimal hypersurface where (2.7) applies and the argument gives via inductive cone reduction that we actually would get a uniform diameter
bound.

\qed

\textbf{Remark} \quad In what follows we can therefore assume that $H$ and $H \setminus \Sigma$ \emph{are connected}, since the subsequent arguments will apply
to
each component.\\

\vspace{1.7cm} \large \textbf {3. Strict positivity of the conformal Laplacian} \\
\normalsize
\bigskip

Now we will construct global conformal deformations $w^{4/n-2} \cdot g$ of $H \setminus \Sigma$ to get metrics with $scal(w^{4/n-2} \cdot g) > 0$. The
interesting feature of this first deformation is its "scaling invariance" close $\Sigma$ which will eventually allow us to induce corresponding solutions on
tangent cones and use them to gain control over the limiting behavior of $w^{4/n-2} \cdot g$ on $H$ near $\Sigma$.\\  For this we will first note that even
weighted forms of the
conformal Laplacian $L$ still have a positive first eigenvalue. \\
This is a property reminiscent of strict stability (cf. [CHS], [Sm]) which would however neither be valid for $H$ in general nor sufficiently versatile since we
want to compare corresponding eigenvalues on $H$ and tangent cones but using the distance to the respective singular sets would not yield steady transitions
since these
sets and therefore the distances change abruptly when passing from even a highly scaled hypersurface to a tangent cone.\\
Instead we will consider versions of (singular) weighted conformal Laplacians $| A | ^{-2} \cdot (- \triangle + \gamma \scal_H)$ where $| A | ^{-1}$ turns out to be the
adequate distance measurement.\\

To see the problem recall that $H \subset M$ is an area minimizer in its homology class and hence stable. Therefore we get from (A2) in the introduction and the
fact that $\scal_M
>0$ that
\[ \inf_{f \not\equiv 0, \supp f \subset H \setminus \Sigma} \frac {\int_{H \setminus \Sigma} |\nabla f|^2 + \frac{n-2}{4 (n-1)} \scal_H f^2} {\int_{H \setminus \Sigma} f^2} > 0.
\]
Beside the problem that $H \setminus \Sigma$ is non-complete such that the spectral properties (and in particular the existence etc. of eigenfunctions) of $L$ will
significantly differ from the closed case, the fact that this \emph{infimum} is positive is just not enough to understand the bending effects of our
construction on the way to the doubled (truncated) hypersurface. One has to use the local/infinitesimal geometry near $\Sigma$. But the eigenvalue is {\bf
not} invariant under scaling and in fact under typical rescaling constructions
as used in the definition of tangent cones the first eigenvalue of the "standard" conformal Laplacian converges to zero. \\

Thus we turn to the weighted operator. The first thing to note is that the metric $g$ on $M$ can be perturbed arbitrarily small in any $C^k$-topology (in
particular keeping $Scal
> 0$) such that the area minimizer $H \subset M$ under consideration does not contain any piece of a totally geodesic hypersurface. More precisely, the set
where $|A| \equiv 0$ is a set of $(n-1)$-dimensional measure zero: Namely we can turn $g$ and therefore $H$ into analytic
objects: then $|A|^2$ becomes an analytic function. If $|A|^2$ is identically zero each tangent cone is regular and hence $\Sigma = \emptyset$. \\

Thus we will henceforth assume that $|A|^2$ is \emph{analytic} and the set \emph{$|A|^{-1}(0)$ has $(n-1)$-dimensional measure zero}.\\

Under this assumption we can derive a {\it strict positivity} property of the conformal Laplacian from its (ordinary) stability:\\

 \textbf{Lemma (3.1)}{\itshape \quad
\[ \lambda_0 := \inf_{f \not\equiv 0,smooth, \supp f \subset H \setminus \Sigma}
\frac {\int_{H \setminus \Sigma} |\nabla f|^2 + \frac{n-2}{4 (n-1)} \scal_H f^2} {\int_{H \setminus \Sigma} | A | ^2  \cdot f^2} > 1/4 \] and we can find a
smooth function positive (although not-integrable) function $u_0$ on $H \setminus \Sigma$ with
\[ - \triangle u_0 + \frac{n-2}{4 (n-1)}  \scal_H u_0 = \lambda_0 \cdot  | A | ^2  \cdot u_0\]} \\

\textbf{Proof} \quad The stability inequality $(A2)$ and $scal_M  > 0$ provide us with the following estimate: \\
\[ \int_H | \nabla f |^2 + \frac{n-2}{4 (n-1)} scal_H f^2 d A \ge \]
\[\int_H \frac{n}{2 (n-1)} | \nabla f |^2 + \frac{n- 2}{2 (n-1)} f^2 \left( | A |^2 + scal_M
\right) d A  \ge \int_H  \frac{n- 2}{2 (n-1)} | A |^2 f^2 d A \]\\

which gives the estimate for $\lambda_0$. The weight as well as the underlying space are singular and thus we cannot handle $\lambda_0$ as a first eigenvalue
with a corresponding eigenfunction by standard means. \\

But we can construct such a smooth function $u_0$ approximating the problem by a sequence of regular ones. \\

1. Define an averaged form of $|A|$,  $\pounds^2_\ve (x) := \frac {\ve^2}{dist_H(x,\Sigma)^2} + |A|^2(x)$:  note that in the case of a cone singularity of the
minimal cone $C V^{n-1}$ over a manifold $V^{n-1} \subset S^{n}$ and with $| A_V | = const$ one has $| A_{C V} | (x) = const \cdot dist (x,0)^{-1}$.  $dist_H
(x, \Sigma)$ is Lipschitz but will not be smooth in general, but letting the heat flow slightly deform this function gives a smooth approximation (which can be
made arbitrarily fine when approaching $\Sigma$) with the additional feature of being $C^{l,\beta}$ close to $dist_H (x, \Sigma)$ in those places where it has
this
degree of regularity (cf.[Fr]). \\
In what follows we think of such a fine smooth approximation; in particular when we speak of level sets
$\pounds_{\ve}^{-1}(d)$ which therefore can (generically) be assumed to be \emph{smooth}.\\\\

2. Choose any exhausting sequence of open subsets of $H$ with compact closure $K_m$ in $H$ such that $\p K_m$ is smooth, $K_m \subset K_{m+1}$,
$\bigcup_{m=1}^\infty K_m = H \setminus \Sigma$. For any $\ve >0$ we find a unique first Dirichlet eigenfunction $u_{m,\ve}$ satisfying
\[ - \triangle u_{m,\ve} + \frac{n-2}{4 (n-1)}\scal_H u_{m,\ve} = \lambda_{m,\ve} \cdot \pounds_\ve^2 \cdot u_{m,\ve}, \quad \lambda_{m,\ve} > 0 \]
with $u_{m,\ve} >0$ on $\op{int}K_m$, $u_{m,\ve} \equiv 0$ on $\p K_m$ and $\int_{B_0} \pounds_\ve^2 \cdot u_{m,\ve}^2 = 1$ for a fixed ball $B_0 \subset H
\setminus \Sigma$.\\

 Since the function space grows for increasing $m$, the eigenvalue $\lambda_{m,\ve}$ decreases monotonically as $m \to \infty$ and hence there is a unique
limit $\lambda_{\infty,\ve} = \lim_{m\to\infty} \lambda_{m,\ve} \ge 0$. Also note that $\lambda_{\infty,\ve} = \lambda_\ve$
\[ \lambda_\ve := \inf_{f \not\equiv 0,smooth, \supp f \subset H \setminus \Sigma} \frac {\int_{H \setminus \Sigma} |\nabla f|^2 + \frac{n-2}{4 (n-1)} \scal_H f^2}
{\int_{H \setminus \Sigma} \pounds_\ve^2  \cdot f^2}  \] since for any function $f$ with compact support in $H \setminus \Sigma$ we eventually have $\supp f
\subset K_m$ for sufficiently large $m$.
\\

{\bf Claim} \quad There is a subsequence of $(u_{m,\ve})_m$ that converges in $C^k$ (for any $k$) to a (not necessarily integrable) limit function $u_\ve >0$ on
$H \setminus \Sigma$ satisfying
\[ - \triangle u_\ve + \frac{n-2}{4 (n-1)} \scal_H \cdot u_\ve = \lambda_\ve \cdot \pounds_\ve^2 \cdot u_\ve. \]
(Note that, unlike $\lambda_\ve$, this limit function {\it may} depend on the
choice of $K_m$.) \\

{\bf Proof} \quad This is a standard application of elliptic estimates and Harnack inequalities. Since such arguments will appear several times later on and the
smoothed weight $\pounds^2_\ve (x)$ might appear unusual, we carry them out
in some detail here. \\

First of all, notice that $\lambda_{m,\ve} \to \lambda_\ve \ge 0$ implies that there exists $c_1 >0$ such that $0 \le \lambda_{m,\ve} \le c_1$ for all $m$.
Thus, on every ball $B$ with compact closure in $H \setminus \Sigma$ the equations
\[ - \triangle u_{m,\ve} + \left( \scal_H - \lambda_{m,\ve} \cdot \pounds_\ve^2 \right) \cdot u_{m,\ve} = 0 \]
have uniformly (in $m$) bounded coefficients. Therefore, we get uniform constants in the interior elliptic estimates
\[ | u_{m,\ve} |_{C^l(B')} \le c_l(B,B') \cdot | u_{m,\ve} |_{L^2(B)} \]
(the $L^2$- and $L_\ve^2$-norms are locally equivalent) and the Harnack inequalities
\[ \sup_{B'} u_{m,\ve} \le \bar c(B,B') \cdot \inf_{B'} u_{m,\ve} \]
for all balls $B' \subset \! \subset B \subset \! \subset H \setminus \Sigma$. \\

Thus, on $B_0$, the $L_\ve^2$-bound = 1 and Harnack's inequality imply upper and lower bounds
\[ c_2'(B_0) > \sup_{B_0} u_{m,\ve} \ge \inf_{B_0} u_{m,\ve} > c_2(B_0) >0 \]
and therefore on a slightly larger ball $B_0' \supset \! \supset B_0$
\[ \sup_{B_0'} u_{m,\ve} \le  c_3 \cdot \inf_{B_0'} u_{m,\ve} \le c_3 \cdot c_2'(B_0), \]
i.e., there is a uniform $L^2$-bound on $B_0'$ and thus a $C^l$-bound on $B_0$ and we may assume that $u_{m,\ve}$ converges in $C^l$ on $B_0$. The limit satisfies
$u_\ve \ge c_0(B_0)
>0$ and the equation
\[ - \triangle u_\ve + \left( \scal_H - \lambda_\ve \cdot \pounds_\ve^2 \right) \cdot u_\ve = 0. \]
\qed

Now is $H \setminus \Sigma$ is connected and hence for any point $x \in H \setminus \Sigma$ outside $B_0$ we can argue by choosing a smooth path $\gamma:[0,1]
\to H \setminus \Sigma$, $\gamma(0) \in B_0$, $\gamma(1)=x$ covered by finitely many overlapping balls $B_1,\dots,B_k$ in order to get $L^2$-estimates: say $B_0
\cap B_1 \neq \emptyset$; then
\[ \tilde c^{-1} \cdot \sup_{B_1} u_{m,\ve} \le \inf_{B_1} u_{m,\ve} \le \inf_{B_0 \cap B_1} u_{m,\ve} \le \]
\[ \sup_{B_0 \cap B_1} u_{m,\ve} \le \sup_{B_0} u_{m,\ve} \le c_3 \cdot \inf_{B_0} u_{m,\ve} \le c_3 \cdot c_2'(B_0). \]\\
Arguing as for $B_0$ we get a further positively lower and upper bounded  subsequence converging on $B_0 \cup B_1$ and, proceeding by induction, a subsequence
converging in $C^k$ to a limit function $u_\ve >0$ on all of $H \setminus \Sigma$. \\
Next we observe that $\lambda_\ve \ra \lambda_0$ for $\ve \ra 0$ and choosing suitable multiples we may assume that $\int_{B_0} u_{\ve}^2 =1$ for every $\ve >
0$. Thus we can argue similarily as before and find a $C^k$-converging sequence $ u_{\ve_i}$ for some sequence $\ve_i \ra 0$, $i \ra \infty$ with smooth limit
$u_0 > 0$ on $H \setminus \Sigma$ satisfying
\[ - \triangle u_0 + \frac{n-2}{4 (n-1)} \cdot scal_H  \cdot  u_0 = \lambda_0 \cdot  | A | ^2  \cdot u_0\] \qed

In a sharp contrast to closed manifolds the previous argument also provides us with an important extension which we use
henceforth throughout many analytic arguments: \\

 \textbf{Corollary (3.2)}{\itshape \quad We can find smooth positive functions $u_\lambda$ on $H \setminus \Sigma$ with
\[ - \triangle u_\lambda + \frac{n-2}{4 (n-1)}  \scal_H u_\lambda = \lambda \cdot  | A | ^2  \cdot u_\lambda\]
for \emph{any} $\lambda <\lambda_0 $. }\\

\textbf{Proof} \quad  We can  decrease the scalar curvature (as described in [L4] ) in each step close to $\p K_m$ for the exhausting sequence $K_m$ such that
for any $\ve>0$ we find for any given $\zeta > 0$ a unique first Dirichlet eigenfunction $u^\zeta_{m,\ve}$ satisfying
\[ - \triangle u^\zeta_{m,\ve} + \frac{n-2}{4 (n-1)} \scal_H u^\zeta_{m,\ve} = (\lambda_{m,\ve}-\zeta  )\cdot \pounds_\ve^2 \cdot u^\zeta_{m,\ve} \]
This leads directly to the first claim.  \qed

 \textbf{Remark} \quad There are two essential points where (3.2) comes into play \\

 1. Since almost all the proofs are by induction we are led to consider  the equation
\[ - \triangle w + \frac{n-2}{4 (n-1)}  \scal w = \lambda_0 \cdot  | A | ^2  \cdot w\] also on a lower dimensional
space, although for the geometric counterpart $n-2/4 (n-1)$ had to be substituted iteratively for $n-3/4 (n-2)$ etc.\\
 While the existence of
solutions is induced from those of the top dimensional equation we had to prove all the properties (iterating the induction) for all these dimensionally shifted
equations.\\
But for $n > 3$ we have $1/6 \le n-2/4 (n-1) \le 1/4$  and since we are seriously concerned only with the case where $ scal = - | A |^2$  we can subsume these
equations under the same label of the dimensionally correct weighted conformal Laplacian but with  a \emph{larger} eigenvalue.\\

2. (3.2)  and more technically refined forms in [L3] will be used to obtain supersolutions of $ - \triangle u_{\lambda_a} + \frac{n-2}{4 (n-1)} \scal_H
u_{\lambda_a} = \lambda_a \cdot | A | ^2 \cdot u_{\lambda_a}$ with transparent properties and here we want to \emph{raise} $\lambda$. \\

 Thus, (for 1. notice  $0 \le n-2/4 (n-1) -
n-2-k/4 (n-1-k) < 1/12$ for $n-k > 3$),
 we will henceforth use and discuss solutions not for the (formal) first eigenvalue $\lambda_0$
but for $\lambda \in [\lambda_0/2, \lambda_0]$  and formulate the top dimensional equations for \[\lambda^0 := \lambda_0/2
> 1/8\] and we choose one fixed smooth
 $ u_{\lambda^0 } > 0$ with $ - \triangle u_{\lambda^0 } + \frac{n-2}{4 (n-1)}  \scal_H u_{\lambda^0 }
= {\lambda^0 } \cdot  | A | ^2 \cdot u_{\lambda^0 }$  \\

\vspace{1.7cm}

\large \textbf {4. Distinguished Eigenfunctions near $\Sigma$} \\
\normalsize

In this section we will modify the previously constructed $ u_{\lambda^0 } > 0$ near the singular set using a Perron-type construction. The new function will
have the particular property of being the {\it smallest} positive eigenfunction with respect to its boundary data and we will find that the minimal solutions
 descend to such minimal solutions on cones. \\

We start with the construction on $H$. The first point to note is that one cannot apply the standard Perron strategy: in general there would not be a solution
to our problem (note that the sign of the linear term is just the converse of the classical case (cf. [GT], p. 103)) and our domain is not complete. \\
But in our case where we have positive boundary data and at least one positive (super-)solution we can adjust the argument to get a minimal positive solution.\\

To begin with, we fix a smoothly bounded neighborhood $V \subset U_\delta(\Sigma)$ of the singular set $\Sigma$ of $H$ within a $\delta$-distance tube
$U_\delta(\Sigma)$. Choosing $\delta \ll 1$ means that $scal_H|_V$ is almost negative: since $ scal_H =  scal_M -
 2 Ric_M (\nu, \nu) - | A |^2 $ the scalar curvature is uniformly upper bounded everywhere and since $\delta \ll 1$ means that in (eventually) most places $|A| \gg 1$
we can scale the whole setting keeping $scal \ll -1$ in most places while $scal_H|_V \ll 1$ everywhere. This can readily be turned into a quantitative statement
using tangent cones where $scal \le 0$ and the zero set is lower dimensional. On $V$ we want to find the smallest solution $u
> 0 $ of the equation (LO)
\[    \triangle u + (\lambda^0 |A|^2- \frac{n-2}{4 (n-1)} \cdot scal_H ) u = 0 \mbox{ with } u \equiv u_{\lambda^0 } \mbox{ on } \p V. \]
For small\emph{ regular }regions the corresponding problem is well-behaved: we know that the first Dirichlet eigenvalue of $\triangle + (\lambda^0 |A|^2-
\frac{n-2}{4 (n-1)} \cdot scal_H )$
will be rather large and this allows us to deduce:\\

{\bf Lemma (4.1)} {\itshape \quad For any $p \in H \setminus \Sigma$ there is a small $\Lambda(p)>0$ such that for any $R \in (0, \Lambda(p))$ the problem
\[ \triangle u + (\lambda^0 |A|^2 - \frac{n-2}{4 (n-1)} \cdot scal_H ) u = 0 \mbox{ on } B_R(p)  \mbox{ and }  u = \varphi \mbox{ on } \p B_R(p) \]

has a unique solution for every (continuous) function $\varphi : \p B_R(p) \to \R$. }\\

{\bf Proof} \quad Choose $R>0$ small enough so that the first Dirichlet eigenvalue $\mu$ of the Laplacian on $B_R(p)$ satisfies $\mu > (\lambda^0|A|^2 -
\frac{n-2}{4 (n-1)} \cdot scal_H )|_{B_R(p)}$. Then the problem
\[ \triangle u + (\lambda^0 |A|^2 - \frac{n-2}{4 (n-1)} \cdot scal_H ) u = 0 \mbox{ on } B_R  \mbox{ and }   u = 0 \mbox{ on } \p B_R(p) \]

has only the trivial solution: for an assumed non-trivial solution $v$ we would obtain
\begin{eqnarray*}
0 &=& \int_{B_R} v \triangle v + (\lambda^0 |A|^2 - \frac{n-2}{4 (n-1)} \cdot scal_H ) v^2  \\
&=& - \int_{B_R} |\nabla v|^2 + \int_{B_R} (\lambda^0 |A|^2 - \frac{n-2}{4 (n-1)} \cdot scal_H ) v^2
\end{eqnarray*}

and therefore
\[ \int |\nabla v|^2 / \int v^2 \le c < \mu \]
contradicting the fact that  $\mu$ gave the minimal value for such quotients.\\
 Thus, using a Fredholm alternative for such elliptic operators ([GT], p.107), the claim follows. \qed

This is used in the study of Perron families on $V$. As usual we
call a function $v:V \to \R$ \emph{supersolution} of $\triangle u
+ (\lambda^0 |A|^2- \scal) u = 0$ if, for any ball $B \subset V$
and any solution $u$  on $B$ with $u|_{\p B} \le v|_{\p B}$, it
follows that $u|_B \le v|_B$.\\
In order to ensure that we have got a sufficiently rich class of supersolutions we first notice that the \emph{minimum} of two supersolutions is also a
supersolution.
\\Now the point is that another operation (a local
upgrading of a super- to an actual positive solution) within this class may at best be valid for small balls.

Thus we first prove the validity and then that this is already
sufficient.\\
The previous lemma allows us: we define the {\it lift} $\bar u$ on
$B_R(p) \subset V$, $R \in (0, \Lambda(p))$ of a supersolution
$u:V\to\R$, as follows: on $B_R(p)$ we let $\bar u$ be the
solution of $ \triangle \bar u + (\lambda^0 |A|^2- \frac{n-2}{4
(n-1)} \cdot scal_H ) \bar u = 0$ with $\bar u|_{\p B_R(p)} =
u|_{\p B_R(p)}$ and $\bar u = u$ on
$V \setminus B_R(p)$. \\

{\bf Lemma (4.2)} {\itshape \quad The lift $\bar u$ of a positive supersolution $u$  is still a positive supersolution. }\\

{\bf Proof} \quad Let $B \subset V$ any ball and consider a solution $h$ of equation (LO) with $h \le \bar u$ on $\p B$.\\. Since $u$ was a supersolution, $\bar
u \le u$ on $B_R(p)$ (and equal outside) and thus $h \le u$ on $\p B$. hence $h \le u$ on $B$ and $h \le \bar u$ on $B \setminus B_R(p)$ and thus $h \le \bar u$
on $\p (B \cap B_R(p))$. But the argument of (4.1) also covers the unique solvability on $B \cap B_R(p)$ and also we will see now that we can handle this case
as soon as we understood $B_R(p)$. Thus eventually we observe that we can reduce the
problem to the following one:\\
 Given a solution $u$ of (LO)  on $B_R$ with $u>0$ on $\p B_R$ we have to show that $u>0$ on all of
 $B_R$.\\
 Consider the family of equations (LOT)
\[ \triangle u + t (\lambda^0 |A|^2 - \frac{n-2}{4 (n-1)} \cdot scal_H ) u = 0. \]
For the moment, consider the situation where $\scal \le 0$ on all of $B_R$. The argument in (4.1) shows again that if $\Lambda(p)$ is chosen small enough we may
assume that for $t \in [0,1]$ (LOT) always has a unique solution
 $u_t$ with $u_t|_{\p B_R} = \varphi
> 0$ and from this we may infer that $u_t$ depends continuously on $t$.\\

 For $u_0$ the claim holds by the minimum principle for harmonic functions. So let
us assume that there are $t \in [0,1]$ and $x_0 \in B_R$ such that $u_t(x_0) < \min \varphi$.\\

Hence we may assume that $0 < u_t(x_0) = \min_{x\in B_R} u_t < \min \varphi$, contradicting the maximum principle, since for
$\lambda^0|A|^2 - \frac{n-2}{4(n-1)}\scal_H \ge 0$, the solution $u_t$ cannot have a positive minimum. \\
In the general situation, i.e. admitting $\scal > 0$ on $B_R$, we argue as follows: Let $u$ be the solution of
\[ \triangle u + (\lambda^0 |A|^2 - \frac{n-2}{4 (n-1)} \cdot \scal^- - \frac{n-2}{4 (n-1)} \cdot \scal^+) = 0 \]
with $u|_{\p B_R} = \varphi$ where $\scal^- = \min (0,\scal)$ and $\scal^+ = \max (0,\scal)$. Let $u_0$ be the solution of this equation with $u_0|_{\p B_r}
\equiv \ve$, where $0 < \ve < \min \varphi$. By uniqueness of solutions on $B_R$ (and any subdomain of it), we obtain $u \ge u_0$ on $B_R$. Moreover, by the
continuous dependence of $u_0$ on the "perturbation" $\scal^+$, we see that for sufficiently small $\scal^+$ both $u_0$ and $u$ remain positive. Now, since the
considered equation is invariant under rescaling, we can, by choosing $V$ small enough, ensure that $\scal^+$ gets arbitrarily small: in $scal_H =  scal_M - 2
Ric_M (\nu, \nu) - | A |^2 $ both $|scal_M|$ and $ | Ric_M (\nu, \nu)|$ decrease \emph{uniformly} quadratically under scaling. \\ More formally one could use a
covering of $\Sigma$ by small balls as obtained in sec. 3 and note that if we choose $V$ as the union of such balls one has a arbitrarily good smooth
approximation by singular cones
outside some small part close to the singularities of the cones. But singular cones (in $\R^n$ have $scal < 0$ almost everywhere) \qed \\

After these preliminary considerations, we are now ready to apply the Perron method to our equation. \\To this end, let $S = \{ v:V \to R \,|\, v \mbox{
supersolution}, v > 0,
v|_{\p V} \ge u \}$; since at least $u \in S$, it is \emph{non-empty}. \\

{\bf Lemma (4.3)}{\itshape \quad  The function $w(x) = \inf_{v\in S} v(x)$ is positive and satisfies}
\[    \triangle w + (\lambda^0 |A|^2- \frac{n-2}{4 (n-1)} \cdot scal_H ) w = 0 \mbox{ with } w \equiv u_0 \mbox{ on } \p V. \]

{\bf Proof} \quad Obviously $w$ is well defined and non-negative. Let $y$ be an arbitrary point of $V$ where $\scal(y) \leq 0$ and $v_n \in S$ such that $v_n(y)
\to w(y)$. By definition, $v_n > 0$ and taking minima (i.e. replacing $v_n$ by $\min (v_n, v_0)$) we may assume that the sequence $v_n$ is bounded. Now choose
$R \in (0, \Lambda(y))$ and assume after scaling and adjusting $V$ that the scalar curvature is almost negative $scal \ll 1$ (in the sense of the previous
proof) and consider the lift $V_n$ of $v_n$ on $B_R(y)$. By lemma (4.2) we have $V_n \in S$ and therefore $w(y) \le V_n(y) \le v_n (y) \to w(y)$. Moreover, by
standard compactness results, we can assume that $V_n$ converges uniformly on any ball $B_\rho(y)$ ($\rho < R$) to an eigenfunction $v$ on $B_R(y)$. Clearly $v
\ge w$ and $v(y)=w(y)$; we wish to prove that $w=v$ on $B_R(y)$: So assume there exists $z \in B_R(y)$ such that $v(z)
> w(z)$. Choose a function $W \in S$ such that $w(z) \le W(z) < v(z)$ and define $w_k = \min (W, v_k) \in S$ as well as the corresponding lifts
$\Bar w_k$ on $B_R(y)$. As before we can assume that $\Bar w_k$
converges to an eigenfunction $\bar w$ on $B_R(y)$ satisfying $w
\le \Bar w \le v$ with equality holding at the point $y$. Hopf's
maximum principle (cf. remark below) gives a contradiction and we
conclude that
$v=w$. \\
It remains to show that $w$ is nowhere zero. To see this, choose a point $x_0 \in \p V$ with $\scal(x_0) < 0$ and a sufficiently small ball $B_R(x_0)$ on which
we have unique solvability of the eigenvalue equation as well as $\scal|_{B_R(x_0)} < 0$. Let $u_0$ be the solution of the equation with boundary data given by
a smooth function $ \phi \ge 0 $ on $(\p B_R(x_0)\cap int V)  \cup  (B_R(x_0) \cap \p  V)$:
\[    \phi \equiv 0 \mbox{ on  } \p B_R(x_0) \cap int V \mbox{ and }  \phi \equiv u \mbox{ near } x_0. \]
Then again by Hopf's maximum principle $u_0 > 0$  on $ B_R(x_0)\cap int V$ and since $v \ge u_0$ for every $v\in S$ we have $w > 0$ on $ B_R(x_0)\cap int V$ and
joining any point in $V$ by a chain of balls we analogously get $w > 0$ on $V$. \qed

{\bf Remark} \quad In order to avoid confusions: the Hopf's maximum principle applies to general solutions of $\Delta u + g(x) u = 0$ with $g \le 0$, which is
precisely \emph{not} our case.\\ But if u \emph{vanishes} in the point where it is applied one can drop the sign assumption for $g$ (cf. [G], p.34) and still
obtains the critical strict inequality for the \emph{outer normal derivative} $\p u / \p n > 0$  in an extremal point $q$ of the zero set in the sense that the
interior ball condition for the complement is satisfied and thus there is a locally (at least relative to this interior ball) unique maximum in $q$.\\ Here and
later on we merge this with non-negativity information to utilize this key estimate (from the proof of Hopf's maximum principle) also for our equations.
We just refer to it as the Hopf's maximum principle.\\

$u_0$ and $w$ can now be patched together along $\p V$  defining a positive continuous function $g_V$ on $H \setminus \Sigma$ with $g_V|_{H \setminus V} \equiv
u_0$ and $g_V|_V \equiv w$. On both subsets $g_V$ is smooth and the conformal deformation on the respective pieces obviously lead to $scal > 0$-geometries. It is
important to note that we can smooth $g_V$ arbitrarily close to $\p V$, say within a neighborhood $W \supset \p V$ in such a way that conformal deformation via
the smoothed function $g^W_V$ on $H \setminus \Sigma$ still gives positive scalar curvature. Actually we will see that we even gain a little bit of positivity if
there is a real crease.\\

{\bf Lemma (4.4)}{\itshape \quad For any neighborhood $W \supset \p V$  we can find a smoothing $g^W_V > 0$ coinciding with $g_V$ on $H \setminus W$ such that $
-\triangle g^W_V + \frac{n-2}{4 (n-1)}  \scal_H g^W_V  > 0$}. \\

 {\bf Proof:} The restriction of $g_V$ to $\p V$ and to distance sets of $\p V$ (close enough they are still submanifolds) is smooth and we only have to care
 about the normal derivatives (directed towards $\Sigma$) of $u_0$ and $w$ along  $\p V$:
 $\p V$ is compact and we may assume that it is connected and claim if $g_V$ was not smooth (that if $w$ is not just $u_0$) then
 $\frac{\p u_0}{\p n} > \frac{\p w}{\p n}$ on $\p V$. Namely we can simply consider $f = u_0 - w$; this is a non-negative somewhere positive solution of
\[ \triangle f + (\lambda^0 |A|^2 - \frac{n-2}{4 (n-1)} \cdot scal_H ) f = 0 \mbox{ on } V  \mbox{ and }  f = 0 \mbox{ on } \p V. \]
Using the proof  of Hopf's maximum principle (cf. [GT], p.34) we first get $f > 0$ on $intV$ and then $\frac{\p f}{\p n} > 0$ on $\p V$.\\
Multiplying $u_0$ by a constant slightly smaller than 1 we can find also shift the original set  $\p V$ where $u_0 \equiv w$ a bit towards $\Sigma$. That means
we can assume that $w$ can be smoothly extend over $\p V$ as a solution of $ -\triangle g^W_V + \frac{n-2}{4 (n-1)}  \scal_H g^W_V  = 0$
and, near $\p V$:  $\frac{\p (u_0 - w)}{\p n} > \kappa > 0$, in particular $u_0 - w > 0$ in $int V$ and $< 0$ in $H \setminus \bar V$. \\

We will meet such a situation later on again and thus we formulate the actual smoothing procedure as an auxiliary\\

{\bf Lemma (4.5)}{\itshape \quad Let $N^{n-1}$ be a smooth submanifold in an orientable manifold $(F^{n}, g)$ and $f_i > 0, i =1,2$ smooth coinciding on
$N^{n-1}$ with $-\triangle f_i + \frac{n-2}{4 (n-1)}  \scal_F f_i  > 0$ on $(F^{n}, g)$ and such that $\frac{f_1}{\p n} > \frac{f_2}{\p n}$ on $N^{n-1}$ .\\
Then for any neighborhood $E$ of $N^{n-1}$ we find a smooth function $f_E > 0$ with $-\triangle f_E + \frac{n-2}{4 (n-1)} \scal_F f_E  > 0$ and, outside $E$,
$f_E = f_1$ resp. $f_E = f_2$ on that side of $N^{n-1}$ where the resp. $f_i$
is the smaller one.} \\

{\bf Proof} \quad Take Fermi coordinates $x_1,..x_n$ in some point
in $N^{n-1}$ such that $x_1 $ is the (unit speed) coordinate in
normal direction.
Then on $N^{n-1}$: $g_{1k} = g^{1k} = \delta_{1k}$ and $det(g_{ij}) = 1$. \\

 The standard Laplacian written in these local coordinates\\ $\triangle f = \frac{1}{\sqrt{det(g_{ij})}} \sum^{n}_{\mu = 1} \sum^{n}_{\nu = 1} \frac{\p}{\p x_\nu}
(\frac{\p f}{\p x_\mu} g^{\mu \nu}\sqrt{det(g_{ij})})$ gives the following formal shape $(ABC)$ for the conformal Laplacian
\[  -\triangle f + \frac{n-2}{4 (n-1)}  \scal_F f =   \sum^{n}_{\mu = 1} \sum^{n}_{\nu = 1} a_{\nu \mu}(g)
\cdot \frac{\p^2 f}{\p x_\nu \p x_\mu} + \sum^{n}_{\mu = 1}
b_{\mu}(g) \cdot  \frac{\p f}{\p x_\mu } + c(g) \cdot f\] where
the coefficients $a_{\nu \mu}(g), b_{\mu}(g), c(g)$ depend only on
the metric and its derivatives and with $a_{11}(g) = 1$ on
$N^{n-1}$ and thus we may assume that $a_{11}(g)\in [1/2,2]$
and that all the other coefficients are uniformly upper bounded on those neighborhoods we are about to choose.\\

For $\delta \ll 1$ the $\delta$-tube $U_{\delta}(N^{n-1})$ is topologically a product and the distance sets (signed distance depending on the side of $N^{n-1}$)
are parametrized by the coordinate $x_1$.

Now take a function $\chi \in C^{\infty}(\R,\R^{\ge 0})$ with $\chi = 0  \mbox{ on }  \R^{\ge 1}$ and $\chi > 0  \mbox{ on }   \R^{< 1}$ with  $\chi(t) >
0$,$\frac{\p \chi(t)}{\p t } < 0 $, $\frac{\p^2 \chi(t)}{\p t^2 } > 0$.
 Moreover for any given positive $K \gg 1$ we can choose $\chi$ such that on $(-1,1)$:
 \[\frac{\p^2 \chi(t)}{\p t^2 } \ge - K \cdot \frac{\p \chi(t)}{\p t } \mbox{ and }\frac{\p^2 \chi(t)}{\p t^2 }\ge K \cdot \chi(t) \mbox{ and  }\chi(0) = 1\]

With  $\chi^\pm_\delta(x) = \chi(\pm x_1/\delta)$ consider for some fixed smooth function $\eta > 0$ on $N^{n-1}$: $f_2 - \eta \cdot \chi^+_\delta(x)$ resp.
$f_1 - \eta \cdot \chi^-_\delta(x)$. For sufficiently large $K$  we observe from $(ABC)$ that both functions still satisfy $ -\triangle f + \frac{n-2}{4 (n-1)}
\scal_F f > 0$ on $U_{\delta}(N^{n-1}) $: the only term that contains $\p^2 \chi(t)/\p t^2 $ is $\p^2(- a_{11}(g) \cdot \eta \cdot \chi^\pm_\delta(x))/\p x^2_1
$. This dominates all the other additional contributions of $\chi$ since in all other terms only the zeroth and first derivatives enter (linearly) and also the
derivatives of $\eta$ along $N^{n-1}$ remain bounded.

For small $\eta > 0$ we can get $f_2 - \eta \cdot \chi^+_\delta(x) > 0$ resp. $f_1 - \eta \cdot \chi^-_\delta(x) > 0$ on $U_{\delta}(N^{n-1}) $ and $\frac{\p
(f_1-f_2)(x) }{\p x_1 } =  \frac{\p (\eta \cdot(\chi^-_\delta(x)-\chi^+_\delta(x)))}{\p x_1 } $ on $N^{n-1}$. Now $K$ depends on $\eta$ but this loop ends if we
now increase $K$; we just need to multiply the fixed $\eta > 0$ by some small constant and may keep $K$. \\
 Thus the functions $f_2 - \eta \cdot \chi^+_\delta(x) $ resp. $f_1 - \eta \cdot \chi^-_\delta(x) $ on the resp. side of $N^{n-1}$
 obviously fit together forming a $C^1$ function which satisfies the requirements outside $N^{n-1}$. Finally we can use a cut-off construction
 $\phi(f_2 - \eta \cdot \chi^+_\delta(x)) + (1-\phi)(f_1 - \eta \cdot \chi^-_\delta(x)) $ to make this a globally smooth function as in our claim:
close to $0$ the zeroth and first derivatives almost coincide and the only second derivative terms that deviate appear as $\phi \cdot \p^2 (f_2 - \eta \cdot
\chi^+_\delta(x))/\p x^2_1 + (1-\phi) \cdot \p^2 (f_1 - \eta \cdot \chi^-_\delta(x))/\p x^2_1 $ but since both function satisfy the inequality this still holds
for this pointwise linear combination. \qed \\

This also concludes the \emph{proof} of lemma (4.4). In what follows we always assume that $W$ has been chosen narrow enough for the subsequent argument so that
we do usually not need to specify $W$ and we call this
resulting  "modified eigenfunction" (which is of course only a supersolution in the creasing area) $u_{mod}$.\\

One uses the Allard approximation of $H$ by tangent cones acc (2.5) to compare $u_{mod}$ with Perron solutions on tangent cones. This helps to understand
$u_{mod}$ close to $\Sigma$. But it involves a comparison of how these cone solutions evolve from solutions defined on compact subsets on $H$ since these almost
isometric sets form (via Allard) the geometric bridge between these spaces. The delicate point is that the limit processes are a priori \emph{non-uniform} and
hence such a comparison argument (which would correspond to an exchange of the order of taking limits)
does not necessarily pass to the limits but we will gain enough control needed for our purposes.\\

The first step is to use the following technical but versatile generalization of the whole construction :\\
Instead of $H$ we can consider $H \setminus W_{k}(\Sigma)$ where  $W_{k}(\Sigma)$ is a smoothly bounded neighborhood of $\Sigma$ with $W_{k+1}(\Sigma) \subset
W_{k}(\Sigma)$ and $\bigcap_k W_{k}(\Sigma) = \Sigma$ . Performing the same constructions as above on $H \setminus W_{k}(\Sigma)$ we also get a solution $u^{H
\setminus W_{k}(\Sigma)}_{mod}$. This positive solution has the following properties: $u^{H \setminus W_{k}(\Sigma)}_{mod} \le u_{mod}$ (since the space of
admissible functions on $H \setminus W_{k}(\Sigma)$ is larger) and for $k \ra \infty$, $u^{H \setminus W_{k}(\Sigma)}_{mod} \ra u_{mod}$ $C^k$-compactly on $H
\setminus \Sigma$. We call the $u^{H \setminus W_{k}(\Sigma)}_{mod}$ \emph{subspace Perron solutions}. \\

These functions can be used to define and approximate Perron solutions of the equation which (since the ambient space is now the flat $\R^n$) looks like
$\triangle \varphi +(\frac{n-2}{4 (n-1)} + \lambda^0) \cdot |A|^2 \cdot \varphi = 0$ on any tangent cone $C^\eta_p$ in some $p \in \Sigma$:\\

{\bf Proposition (4.6)} \quad \emph{  For any $\delta > 0$ and any triple $R \gg 1 \gg \varrho \gg \xi >0$ we can find a small $\eta_{\delta , R , \varrho, \xi}
> 0$ such that for \textbf{every} $\eta \in (0, \eta_{\delta , R , \varrho, \xi})$ we can assume:
\begin{enumerate}
\item $ |g_\eta|_{C^k} < \delta \mbox{ on } (B_{R}(0) \setminus B_{\varrho}(0)) \setminus V_{\xi}(\sigma) \subset  C^\eta_p$ \item after $L^2$-normalization we
have for $k = K_n \cdot \eta^{-2}$ \[1 - \delta < u^{H \setminus W_{k}(\Sigma)}_{mod} / c(\omega) r^\alpha < 1 + \delta \mbox{ on } (B_{R}(0) \setminus
B_{\varrho}(0)) \setminus V_{\xi}(\sigma)\] for some Perron solution $c(\omega) r^\alpha$ on $C^\eta_p$.
\end{enumerate}
The second statement means that the Perron solution on the cone $c(\omega) r^\alpha$ respects the cone structure, where $(\omega,r) \in C^\eta_p$  are the
spherical and radial coordinates, $c(\omega) > 0$ a smooth function on $\p B_1(0) \cap C^\eta_p$ and there are constants $-\frac{n-2}{2}  < \theta_1(n) <
\theta_2(n) < 0$ such that $\alpha \in (\theta_1(n),\theta_2(n))$ and $\alpha = \alpha_p$ is uniquely determined for every $p \in \Sigma$.\\}

 \emph{ Moreover, this can actually be refined to obtain that these claims still hold for $u_{mod}$ instead of $u^{H \setminus W_{k}(\Sigma)}_{mod}$. }  \\

The reason why we get that Perron solutions on the cones have the special structure $c(\omega) r^\alpha$ is that this minimality can be seen as a kind of
uniqueness
statement and then the symmetry of the cone structure induces such a separation of variables.\\
The estimate  $\alpha  < 0$ means that the cone will become acuter when being deformed with
 $c(\omega) \cdot r^\alpha$ moreover here is a compilation of some other important geometric properties:\\

{\bf Proposition (4.7)}{\itshape \quad   $C$ equipped with the metric $\tilde g := (c(\omega) r^\alpha)^{4/n-2} \cdot g$ is again a \textbf{cone} (although not
embed) with finite distance between $0$ and any other point of $C$:\\ $(C,\tilde g)$ is isometric to any of copy scaled around $0$ and can be reparametrized as
$c(\omega)^{4/n-2}  \cdot g_{\R} + r^2 \cdot g_{\p B_1(0) \cap C}$ and the scalar curvature in a point with new distance $\rho$ to $0$ is equal
to $\frac{4 (n-1)}{2 |\alpha| } \cdot \lambda^0 \cdot c(\omega)^{4 (n-3)/n-2} \cdot a(\omega)^2/\rho^2$. \\
We can conformally deform the metric $\tilde g$ to some other cone metric $\tilde g^\ast$ with
$scal_{\tilde g^\ast}(\omega,\rho) \ge \iota_H/\rho^2$ for some $ \iota_H > 0$ which is independent of the singular cone $C  \in \overline {\cal T}_H$}\\

Actually these deformations can be done in a natural way in the sense that their definition changes continuously on the space of cones in flat norm topology.\\

The inductive way these results are proved also leads to an important extension of the previous statements for the spherical component of a Perron solution
$c(\omega)$ satisfying
\[(CW) \;\; \left(\alpha^2 + (n-2) \alpha \right) \cdot c(\omega) +
\left( \triangle_S + \left( \frac{n-2}{4 (n-1)}+\lambda \right) a(\omega)^2 \right) c(\omega) = 0 \] Note that the first term disappears when passing to tangent
cones of $\p  B_1(0) \cap C$. \emph{$c(\omega)$ is also approximated by Perron solutions for (CW) on cones which again split into radial and spherical
components.} In other words there is an inductive descend via Perron solutions of dimensionally
shifted equations.\\

{\bf Corollary (4.8)}{\itshape \quad $scal(\p B_1(0) \cap C,\tilde g^\ast|_{\p B_1(0) \cap C}) > 0$ near $\sigma$}\\

\textbf{Proof} \quad Since $c(\omega) \ra \infty$ and  $a(\omega) \ra \infty$ almost everywhere near $\sigma$ we notice (in view of the following scalar
curvature redistribution) that $ \left( \triangle_S + \left( \frac{n-2}{4 (n-1)}+\lambda \right) a(\omega)^2 \right) c(\omega) \approx 0$ near $\sigma$. Since
$\frac{n-2}{4 (n-1)} > \frac{n-3}{4 (n-2)}$ we also get $ \left( \triangle_S + \left( \frac{n-3}{4 (n-2)} + \lambda \right) a(\omega)^2 \right) c(\omega)
\approx < 0$ near $\sigma$. But that means that the transformation law (TL) finally gives $scal(\p B_1(0) \cap C,\tilde g^\ast|_{\p B_1(0) \cap C}) > 0$ near
$\sigma$. \qed

{\bf Remark} \quad  The verification of (4.6) and (4.7) needs some involved cone reduction argument isolated from the rest of the present paper and appears as a
part of [L3]. The reader may have noticed the word "assume" in the statement of (4.6): there is also a weaker geometric argument that allows to argue
geometrically that $u_{mod}$ can be bent keeepin $scal > 0$ in a way generalizing (4.4) (usually this is called an \emph{h}-principle) to get a transition in
the limit to
a cone Perron solution.\\

\vspace{1.7cm} \large \textbf {5. Geometry on collections of deformations} \\
\normalsize

Now we prepare our new geometry $u_{mod}^{4/n-2} \cdot g_H$ for a second finer semi-local construction. First note that we may assume that $scal > 0$ and as a
result of the techniques we will explain here one can use (4.7) to deform $u_{mod}^{4/n-2} \cdot g_H$ in such a way that the induced geometry on tangent cones
 also has $scal > 0$ as described in (4.7). This is just a conformal redistribution of scalar curvature and later on we can
  assume that $u_{mod}^{4/n-2} \cdot g_H$ already had these properties to start the more involved barrier set up.\\

The technical clue for these constructions is a covering technique of $\Sigma$ by balls with several features like uniformly good approximation by some cone
geometry and universally upper bounded
intersection number.\\

These balls will be measured with respect to the metrics $u_{mod}^{4/n-2} \cdot g_H$ on $H$ respectively on $C$ equipped with the metric $\tilde g := (c(\omega)
r^\alpha)^{4/n-2} \cdot g$ (actually with $\tilde g^\ast$).\\
Hence we resume the discussion of the distance function from (2.7). Since $u_{mod}$ has a pole in $\Sigma$ the singular set will be stretched to \emph{infinite}
length when \emph{measured intrinsically} (the same applies to $\sigma$ and $c(\omega)$) but since \\
$(r^\alpha)^{2/n-2} \le r^{-\beta_n}$ for $\beta_n < 1$ the distance from any interior point remains finite (using the same network of pathes as in (2.7)) and
therefore the distances on $\Sigma \subset
H$ and $\sigma \subset C$ also remain finite (via shortcuts in the interior of $H$ and $C$).\\

We first discuss the effect of the deformation on the notion of distance balls around the tip of a cone
(since all further conformal deformations are truncated Green's functions
around the tip of some cone):\\
As a result from the inductive proof of (4.6) the pole order of $c(\omega)^{2/n-3}$ is also uniformly smaller than $ -1$: for each tangent cone of $\p B_1(0)
\cap \sigma \setminus \{0\}$  the Perron solutions induced from $c(\omega)$ will have a pole of order $\alpha'$ uniformly bounded within $-\frac{n-3}{2} <
\alpha' < 0$. Thus we can also apply the path network idea of (2.7) also for $\sigma \subset C$ equipped with $\tilde g^\ast$.\\ Quantitatively, we have that
the lower bound for the pole order of the deformations induced inductively on top and all lower dimensional tangent cones is  $\in (-\Theta^-_n, -\Theta^+_n)$
for $0 < \Theta^+_n < \Theta^-_n < 1$. Thus we get for distances in radial direction of a cone: \[k_1(n) \cdot (d_{(C, g_C)}(0,x))^{1-\Theta^+_n} < d_{(C,\tilde
g^\ast)}(0,x) \le k_2(n) \cdot (d_{(C, g_C)}(0,x))^{1-\Theta^-_n}\] for $x = (\omega,r)$ where $k_2(n) \ge k_1(n) > 0$ can be chosen independently of the cone
$C$ within the compact set $SC_n$ which is seen from
a cone reduction since we have a uniform upper bound for the diameter from (2.7). Note that such inequalities do not hold on $H$.\\

Although we will use $u_{mod}^{4/n-2} \cdot g_H$ on $H$ later on we stick to $g_H$ for the following covering argument since this allows us to take over
combinatorial properties of ball coverings in Euclidean spaces from the embedding of $H$ in the smooth $M^{n+1}$. \\Thus we first recall from sec.2 where we
have seen that for each point $p$ in $\Sigma$ there is a individual radius $r(p,Q)
> 0$ such that from this radius on (downwards) every $\rho^{-2}$-rescaled ball $B_\rho$ is approximated by a cone up to some approximation quality $Q$ (to be
specified in terms of the various norms). A covering of balls all with quality $Q$ will allow us to handle a number of constructions in a shell (=difference of
a neighborhood \emph{minus} a smaller neighborhood) of $\Sigma$ of controlled tiny size in a local way that composes from single balls with controllable defects
globally. This will be subject of the later sections. For now we construct coverings with
a number of properties whose further impact is just indicated by some keywords.\\

 We
subdivide $\Sigma$ into parts $\Sigma^a_j$ according to their approachability via tangent cones with the notations of the previous section.
\[\Sigma^a_j = \{ x \in \Sigma \; | \; d(\tau^{-2} \cdot (H \cap B_{\tau}(x)), C_{\tau} \cap B_1(0)) \le 2^{-a} \mbox{\emph{ and} } \]
\[\tau^{-1} \cdot (H \cap B_\tau(x) \setminus V_{ 2^{-2^{a}}}(\sigma) \mbox{ \emph{ can be written as graph of a smooth}  } g_\eta \mbox{\emph{ with }}\]
\[ |g_\eta|_{C^k} <  2^{-a} \mbox{\emph{ on} } (B_{1}(0) \setminus B_{ 2^{-a}}(0)) \setminus V_{ 2^{-2^{a}}}(\sigma)\mbox{ \emph{ for any } } \tau \le 2^{-j} \]
\[\mbox{ \emph{for some suitable tangent cone} }
C_{\tau} \mbox{ of  } H \mbox{  in  } x\} \]
\medskip

From (2.4) and (2.5) we observe that for any $a > 0$ we have \[\Sigma^a_j \subset \Sigma^a_{j+1} \; \mbox{ and } \; \bigcup_j \Sigma^a_j = \Sigma.\] We consider
the difference sets $\Delta\Sigma^a_j = \Sigma^a_{j} \setminus \Sigma^a_{j-1}$, $\Delta\Sigma^a_1 = \Sigma^a_1$.
(From these descriptive definitions one can check that these sets are measurable.) \\

Since the geometric properties used in our covering arguments are local we can assume that the \emph{ambient} manifold $M^{n+1}$ had been scaled such it looks
uniformly nearly flat in
$C^k$-topology at a local level i.e. in what follows we may consider $(M^{n+1},g)$ of being $\R^{n+1}$ with its Euclidean metric.\\
Now we cover the $\Delta\Sigma^a_j$ by distance balls $B_{\varrho_{ij}}(p_{ij})$ whose radius is measured with respect to the ambient space. For our purpose we
can also consider them (up to a uniformly negligible error) as intrinsic distance balls in $H$:  in view of (2.7) this will become clear only when we carry out
the deformations in sec.6 and 7. The point will be that we will be able to ignore certain tiny neighborhoods of $\Sigma$ and outside this set the definition of
$\Delta\Sigma^a_j$ allows us
to identify $H$ with some tangent cone centered in $p_{ij}$ and it will only be here that these distances are of interest.\\

Now we define a covering ${\cal{B}}(a,j,\delta)$ of $\Delta\Sigma^a_j$ by $n+1$-dimensional balls $B_{\varrho_{ij}}(p_{ij}) \subset M^{n+1}$, $p_{ij} \in
\Delta\Sigma^a_j, i \in I_{\delta, j}$ such that $\varrho_{ij} \le \epsilon (p_{ij}) $ for a size parameter $\epsilon (p_{ij}) \in (0,1)$ which depends
(discontinuously) on the base point $p_{ij}$: at this stage we could use $  \epsilon (p_{ij}) = 2^{-j} $ but we choose it later (smaller)
when we apply (4.6) in order to to view $u_{mod}$ as a Perron solution $c(\omega) r^\alpha$ on a tangent cone $C^\eta_p$. \\

For any $z \in B_{\varrho_{ij}}(p_{ij})$ we have $B_{\varrho_{ij}}(p_{ij}) \subset B_{2\varrho_{ij}}(z)$ and we consider the two coverings  \[ F^\ast_a =
\bigcup_{ i \in I_{\delta, j} , j \ge 1 } \{  B_{\varrho_{ij}}(p_{ij})\}\]
\[ F_a = \bigcup_{i \in I_{\delta, j} , j \ge 1 } \{ B_{2\varrho_{ij}}(z) \; | \quad z \in B_{\varrho_{ij}}(p_{ij})\}\]\\
Notice that each section $\Gamma \in \prod_{ i \in I_{\delta, j} , j \ge 1 } \{ B_{2\varrho_{ij}}(z) \; | \; z \in B_{\varrho_{ij}}(p_{ij})\}$ can be considered
as a covering of $\Sigma$.\\

Since the indices $i$ and $j$ of all the radii $\varrho_{ij}$ form a countable set we can assume (by small perturbations) that the various $\varrho_{ij}$ are
pairwise unequal i.e. $\varrho_{ij} =
\varrho_{kl}$ iff $i=k$ and $j=l$.\\
Moreover we can assume such that for each $\varrho_{ij}$ there is a \emph{largest} $\varrho_{i^\ast j^\ast}$ which is smaller than $\varrho_{ij}$ (and in our
case we may also assume there is a largest radius $\varrho_{i^0j^0}$ under all $\varrho_{ij}$) i.e. these radii form a well-ordered set (with \emph{order}
converse to the size of
the radii)  amenable to transfinite induction.\\

Now a Besicovitch style argument cf.([F],Ch.2.8) gives us subcollections of  $F_a$ which still form a covering of $\Sigma$:\\

{\bf Proposition (5.1)} \quad \emph{There are $c(n)$ disjoint families ${\cal{G}}_l \subset F_a$ , $1 \le l \le c(n)$ of balls $B_{2\varrho_{ij}}(z_{ij})$,
$z_{ij} \in B_{\varrho_{ij}}(p_{ij})$ such that for any two balls $B_{2 \varrho_{ij}}(z_{ij}), B_{2 \varrho_{kl}}(z_{kl})$ within one family ${\cal{G}}_{l_0}$,
$B_{10 \cdot \varrho_{ij}}(z_{ij}), B_{10 \cdot \varrho_{kl}}(z_{kl})$ do not intersect and within different families ${\cal{G}}_{l_1}, {\cal{G}}_{l_2}$ the
balls of radii $\varrho_{ij}$  do not contain the center of other balls: $z_{kl} \notin B_{ 2\varrho_{ij}}(z_{ij}), z_{ij} \notin B_{2 \varrho_{kl}}(z_{kl})$,
$z_{kl} \in {\cal{G}}_{l_1}, z_{ij} \in {\cal{G}}_{l_2}, l_1 \neq l_2$ with $\Sigma \subset \bigcup_{{\cal{G}}_l,
1 \le l \le c(n)} B_{2\varrho_{ij}}(z_{ij}) $.}\\

Note that the radii of \emph{these} balls do \emph{not} correspond to the degree of cone approachability around their midpoints but tracing the constructing
back to the points $p_{ij}
\in I_{\delta, j}$ will provide us with this essential feature as well. \\
\medskip

\textbf{Proof of (5.1)} \quad Since $F_a$ also contains the balls (scaled by two) of the covering $F^\ast_a$, we know that the $\varrho_{ij} \ra 0$ uniformly in
$k$ and since  $F_a$ contain all $2\varrho_{ij}$-balls $ B_{2\varrho_{ij}}(z),\; z \in B_{\varrho_{ij}}(p_{ij}), B_{\varrho_{ij}}(p_{ij}) \in F^\ast_a$  we note
that for any $p \in
\Sigma$ : $\inf \{r\; | \; B_{r}(p) \in F_a \} = 0$.\\

Now we construct the families ${\cal{G}}_l$: identifying the set $F_a$ with the well-ordered set of distinct radii of balls we define a map $f: F_a \ra \Z^{\ge
0}$ whose meaning is that its value, say $l$, is the index of the family ${\cal{G}}_l$ where it will belong to. The index $0$ however means that
this ball is ruled out, i.e. is neither used nor needed for the covering.\\
The definition is by (transfinite) induction and starts with the largest radius $\varrho_{i^0j^0}$: we choose $f(\varrho_{i^0j^0}) := 1$ and assume inductively
$f$ had been defined for all radii $ \varrho_{ij}> \varrho_{kl}$. Then we set (which has to justified below) \quad $f(\varrho_{kl}) := $
\[{\footnotesize \left\{%
\begin{array}{ll}
    0 & \hbox{if } z_{kl} \in \bigcup_{\{\varrho_{ij} | \varrho_{ij} > \varrho_{kl}, f(\varrho_{ij}) >0 \}} B_{2\varrho_{ij}}(z_{ij}) \\
min(\{f(\varrho_{ij}) > 0 \; | \; B_{10 \varrho_{ij}}(z_{ij})\cap B_{10 \varrho_{kl}}(z_{kl}) = \emptyset \} \\ \quad \cup \quad \{max \{f(\varrho_{ij})
\; | \;\varrho_{ij}> \varrho_{kl}\}+1\} ) &  \hbox{ otherwise }\\
\end{array}%
\right.}\]

The second option is non-trivial: recall that in $\R^{n+1}$ there is a constant $M(n)$ such that for any configuration of balls $B_{r_k}(z_k)$ with $\|z_k\|
> r_k > 1$ which intersect $B_1(0)$ in such a way that each of the centers $z_k$ is not contained in any other of the intersecting balls  the number
of these balls $B_{r_k}(z_k)$ is at most $M(n)$.\\
Thus if $\varrho_{kl}$ is not an accumulation point within the subset $f^{-1}(\Z^{\ge1}$) \[ min(\{f(\varrho_{ij}) \; | \; B_{10 \cdot \varrho_{ij}}(z_{ij})\cap
B_{10 \cdot \varrho_{kl}}(z_{kl}) = \emptyset \} \cup \{max \{f(\varrho_{ij}) \; | \;\varrho_{ij}> \varrho_{kl}\}+1\} )\] is well defined. Moreover starting
from $\varrho_{i^0j^0}$ whose value under $f$ is $1$ the sequence of values cannot exceed $M(n) +1$ as long as we do not meet an accumulation point: since
otherwise this implies that
the last ball (which is smaller than its proceeders) meets at least $M(n) +1$ balls such that all of them do not contain the center of an other ball.\\
Also for an accumulation point within the subset $f^{-1}(\Z^{\ge1})$ we observe if the value exceeds $M(n) +1$ in $\varrho_{kl}$ this implies it has infinitely
many intersections with larger balls not containing the center of an other ball which cannot happen.\\

Thus we have $c(n) = M(n) +1$ families ${\cal{G}}_l := f^{-1}(l) \subset F_a$, $1 \le l \le c(n)$ of balls $B_{2\varrho_{ij}}(z_{ij})$, $z_{ij} \in
B_{2\varrho_{ij}}(p_{ij})$. By definition they have the properties such that for any two balls $B_{2 \varrho_{ij}}(z_{ij}), B_{ 2\varrho_{kl}}(z_{kl})$ within
one family ${\cal{G}}_{l_0}$, $B_{10 \cdot \varrho_{ij}}(z_{ij}), B_{10 \cdot \varrho_{kl}}(z_{kl})$  do not intersect and within different families
${\cal{G}}_{l_1}, {\cal{G}}_{l_2}$ the balls of radii $2\varrho_{ij}$  do not contain the center of other balls in $f^{-1}(\Z^{\ge1})$. Their union forms a
covering: if $q$ was a point not covered we know from  $\inf \{r\; | \; B_{r}(q) \in F_a \} = 0$ that can take one of the balls in $F_a$ around $q$ and observe
that $f$ would give a value $> 0$ since it would not be ruled out just because its center is not contained in one of the other remaining balls. \qed

\medskip

Now we use that the $\varrho_{ij}$ are pairwise unequal. This allows us to identify the center $p_{ij}$ of balls such that $B_{2\varrho_{ij}}(z_{ij}) \in
\bigcup_{1 \le l \le c(n)}{\cal{G}}_l $ and we immediately get\\

 {\bf Corollary (5.2)} \quad \emph{There are $c(n)$ disjoint families ${\cal{F}}_l \subset F_a$ , $1 \le l \le c(n)$ of balls $B_{
\varrho_{ij}}(p_{ij})$ such that}
\begin{enumerate}
\item for any two balls $B_{ \varrho_{ij}}(p_{ij}), B_{ \varrho_{kl}}(p_{kl})$ within one family ${\cal{F}}_{l_0}$\[B_{6 \cdot \varrho_{ij}}(p_{ij})\cap B_{6
\cdot \varrho_{kl}}(p_{kl}) = \emptyset\] \item within different families ${\cal{F}}_{l_1}, {\cal{F}}_{l_2}$ the balls of radii $\varrho_{ij}$ do not contain
the center of other balls: \[p_{kl} \notin B_{ \varrho_{ij}}(p_{ij}), p_{ij} \notin B_{ \varrho_{kl}}(p_{kl}), \; p_{kl} \in {\cal{F}}_{l_1}, p_{ij} \in
{\cal{F}}_{l_2}, l_1 \neq l_2\]
 \item   $ \; \Sigma \subset \bigcup_{{\cal{F}}_l, 1 \le l \le c(n)} B_{3 \cdot \varrho_{ij}}(p_{ij})$
\end{enumerate}
We finally set  \[{\cal{F}}^a := \bigcup_{1 \le l \le c(n)}{\cal{F}}_l\]

\medskip

\vspace{1.7cm} \large \textbf {6. Surgery in Hausdorff-codimension $> 2$} \\
\normalsize

The upshot of our construction so far is that we can conformally deform $H$ to $scal > 0$ in such a way that pointwise the geometry near $\Sigma$ can be
regarded as a cone geometry with $scal >0$.\\ The next series of deformations (located close to $\Sigma$) bends the geometry in a way that the combination of
all these deformations forms a barrier for $n-1$-dimensional area minimizers within $H$ preventing them from reaching any point of $\Sigma$.
Each of these deformations looks like a Green's function and we place them around points in stratified sets which arise as natural approximations of $\Sigma$.\\

We use this in combination with obstacles placed outside an actually small neighborhood of $\Sigma$ such that the interesting area minimizer (with obstacles)
stays close to $\Sigma$. Otherwise, even any local area minimizer in the homology class could be just a point (since boundaries of neighborhoods of $\Sigma$ are
clearly null-cobordant). Moreover, compressions of minimal hypersurfaces close to
$\Sigma$ will allow us to localize the proof of tightness of barriers.\\

Thus we formalize the notion of a barrier in the following way: Let $G^n \subset F^{n+1}$ be differentiable manifolds, $F^{n+1}$ closed, $G^n$ properly embedded
but incomplete, ${\cal D}= \overline{G} \setminus G$. Only now we also equip $G^n$
with a Riemannian metric (a priori not induced from a metric on $F^{n+1}$)\\

Formally, take two (for now) piecewise smooth compact and cobordant but not necessarily connected submanifolds $M_1^{m}, M_2^{m}$ and the cobordism $W^{m+1}$
equipped with some Riemannian
metric.\\

{\bf Definition (6.1)}{\itshape \quad An area minimizing current ${\cal{T}}$ in $W^{m+1}$ homologous to $M_1^{m}$ (and thus to $M_2^{m}$) is called an area
minimizer
with obstacles $M_1^{m}$ and $ M_2^{m}$.}\\

In most applications one the two obstacles corresponds just to a compactness condition and is never really touched by the support of ${\cal{T}}$ and thus we
will only refer to the effective obstacle as \emph{the} obstacle.\\

 {\bf Definition (6.2)} \quad \emph{
 We call a pair of neighborhoods $(Y,Z)$, $Y \subset \overline{Y} \subset intZ \subset G$ of ${\cal D}$ a \emph{\textbf{barrier}} of ${\cal D}$ if for any
 area minimizer $W ^{n}\subset Z \setminus Y$ with obstacle $\p Y \cup \p Z$ within the same homology class as $\p Z$ the support of $W ^{n}$ does not reach $\p Y$, i.e.
$W ^{n} \cap \overline{Y} = \emptyset$ .}\\

 If $\Sigma$ is a submanifold \emph{and} if there is a neighborhood of  $\Sigma$ in $H$ that looks like a product, then the almost (= after scalings)
Euclidean geometry gives an obvious clue of how to find coverings by deformations of Green's functions type comprising a tight barrier: choose the positions and
coefficients such that in the limiting case of infinitely many superposed deformations one gets the Green's function along $\Sigma$.\\ We observe that the
classical \emph{codim 3 surgery} techniques for $scal> 0$ in the way described in [SY4] appears as a continuous limit of this technique in the special case
where $\Sigma$ is a non-singular submanifold and there is a neighborhood isometric to a product cone. However, from the general viewpoint, the approaches in
[GL1] and [SY4] contain a not really compelling coupling with the simultaneous construction of
 a totally geodesic boundary.   \\

In general there will be configurations of deformations which partially annihilate their deflecting effect: if $\Sigma$ is something like a fractal with
iterated zigzag lines it is easy to find arrangements of deformations where the deflecting effect of one deformation pushes a
minimal hypersurfaces into $\Sigma$ on the opposite side even when we choose arbitrarily small radii.\\

Therefore we use a bit more information about $\Sigma$ than just its codimension: the meta-relation $\Sigma \prec \sigma$ provides us with such a piece of extra
structure that guides us to form more robust local barriers and assemble them to global barriers via the following three level construction: \\

1. We first define single truncated Green's functions (i.e. outside a ball they are extended by 1) on cones \emph{after} having carrying out the deformation
using the Perron solution and estimate the negative impact for the positive scalar curvature in the cut-off region. We call the resulting conformal deformation
(but for simplicity also the function) an \textbf{elementary barrier}, since it is obvious that an area minimizer homologous to the distance unit sphere
will stay away from a region close to the tip of the cone.\\

2. A collection of elementary barriers built from truncated Green's functions is used to assemble the \textbf{local barriers} (still on the cone)
which are composed as follows:\\
For a cone  $C$ with singular set $\sigma$ we start with an elementary deformation around $0$ and after that an area minimizer homologous to the distance unit
sphere will stay outside say  $B_{2r}(0) \cap C$.\\ Now we additionally place families of much smaller elementary barriers along an extension of $\sigma \cap
B_2 \setminus B_r$. We use the existence of global barriers (of sec. 7 below) in the lower dimensional case to get a global barrier for the extension of $\sigma
\cap \p B_1$ (including the inductive definition of this extension).   (Recall from (2,1) that $\sigma \cap \p B_1$ can be handled like the singular set of an
area minimizer). Then one extends this configuration scheme for $\sigma \cap \p B_1$ in cone direction.
These two steps give a barrier around a $\sigma \cap B_1 $ .\\

3. Now we use the coverings with upper bounded intersection number $c(n)$ of sec.5 to transplant this creation (suitably scaled) from tangent cones to
\emph{well} approximated parts of $H$. \\
 The point is that although we will
not use any information about the relative positions of the local barriers (except for the intersection number $c(n)$) their combined deflecting effect can be
estimated as we can refine the approximation (that it we make $\Sigma$ "thinner relative to $\sigma$) while keeping the barriers fixed.\\

This will give the \textbf{global barrier} for $\Sigma$ augmented by certain (germs of) smooth regions in $H$ where $|A|$
grows faster than quadratically which (following the inductive strategy) could be regarded an approximation of $\Sigma$ by \emph{stratified sets}. \\

\bigskip

Now we start with the \emph{elementary barriers}: Let $p \in \Sigma$ be a singular point and $C_p$ a tangent cone at $p$. On $C_p$ we consider the metric
conformally deformed by $u^{4/n-2}$, $u = c(\omega) \cdot r^\alpha$ the Perron solution on $C_p$. And further deformed to a $scal > 0$ cone metric $\tilde
g^\ast$ with $scal_{\tilde g^\ast}(\omega,\varrho) \ge \iota_H/\varrho^2$ for some $ \iota_H > 0$ as in ...We choose the radial distance $\rho$ from $0$ on
$C_p$ measured with respect to the metric
$\tilde g^\ast$  and  consider the "Euclidean" Green's function $\varphi(\rho) = \frac{1}{\rho^{n-2}} $:\\

{\bf Lemma  (6.3)}{\itshape \quad We still have: $\triangle \varphi = 0$ for the Laplacians computed with respect to the metric $\tilde g^\ast$.}\\

\textbf{Proof} \quad We check that also for the abstracted cone geometry $\tilde g^\ast = \tilde c(\omega)^{4/n-2} \cdot g_{\R} + r^2 \cdot g_{\p B_1(0) \cap
C}$ the Laplacian  $\triangle f$ has again the form  $\frac{\p^2 f}{\p \rho^2} + \frac{n-1}{r} \frac{\p f}{\p \rho} + \frac{1}{\rho^2} \triangle_S f$ with $S =
\p B_1(0)$ measured in $\tilde g^\ast$ (here $\triangle_S f$ means the $(n-1)$-dimensional Laplacian of $f(\rho \cdot x)$ calculated on $S = \p B_1(0)$):\\

For $\triangle f = \frac{1}{\sqrt{det(g_{ij})}} \sum^{n}_{\mu = 1} \sum^{n}_{\nu = 1} \frac{\p}{\p x_\nu} (\frac{\p f}{\p x_\mu} g^{\mu \nu}\sqrt{det(g_{ij})})
$ we choose $x_1 = \rho$ and $x_2,.. x_ {n}$ geodesic coordinates for some point $y \in S$  (scaled by $\rho$ when we consider the corresponding point in $\p
B_\rho(0)$ when the
metric is scaled by $\rho^2$).\\
Then we have in the point $\rho \cdot y$: $\sqrt{det(g_{ij},\{i \ge 2, j \ge 2\})} = \sqrt{det(g_{ij})} = \rho^{n-1}$, $ g^{1 k} =  g^{k 1} = 0$ and $\p
g_{ij}/\p x_k = 0$ for $k
> 1$, $g^{11} = g_{11} = 1$  thus
{\small \begin{eqnarray*}  \triangle f &= & \p^2 f/\p \rho^2 + \frac{n-1}{\rho} \cdot \p f/\p \rho + \frac{1}{\sqrt{det(g_{ij})}} \sum^{n}_{\mu = 2}
\sum^{n}_{\nu = 2} \frac{\p}{\p
x_\nu} (\frac{\p f}{\p x_\mu} g^{\mu \nu}\sqrt{det(g_{ij})}) \\
&=& \p^2 f/\p \rho^2 + \frac{n-1}{\rho} \cdot \p f/\p \rho +\frac{1}{\rho^2} \triangle_S f
\end{eqnarray*}}
Hence, $\triangle \varphi = \p^2 \varphi/\p \rho^2 + \frac{n-1}{\rho} \cdot \p \varphi/\p \rho + \frac{1}{\rho^2} \triangle_S \varphi = (n-1)(n-2)\rho^n +
(-(n-1)(n-2)\rho^n) + 0 = 0$.\qed

This function can now be truncated (this allows us to use them on $H$ locally). Since we will combine such functions by multiplication we start with
$\varphi_\mu(x) = \varphi_\mu(\rho)= \mu \cdot \frac{1}{\rho^{n-2}}  \mbox{, with } d(x,0) = \rho,  \mbox{ on } C_p$ and define another radially symmetric
smooth $\phi_{\mu} \ge 1 $ constant $\equiv 1$ outside some ball by a uniform cut-off:
\[\phi_\mu(x) = \left\{%
\begin{array}{ll}
    \varphi_\mu(\rho) & \hbox{on }  B_1(0)  \\
\varphi_\mu(\rho) \cdot \chi(\rho)  & \hbox{on }   C_p \setminus B_1(0)\\
\end{array}%
\right.\]

$\chi$ is defined as follows: choose some fixed function $\chi \in  C^\infty(\R,[0,1])$ on $\R$ with  $\chi \equiv 1$ on $\R^{\le 1}$, $\chi \equiv 0$ on
$\R^{\ge 2}$.\\

For small $\mu$ deforming by $\phi^+_\mu .= \phi_\mu + 1$ will not diminish the positivity of $scal_{\tilde g^\ast}(\omega,\rho) \ge \iota_H/\rho^2$\\

{\bf Lemma (6.4)}{\itshape \quad There is an $\mu_H > 0$ such that for any singular minimal cone $C  \in \overline {\cal T}_H$ equipped with $g_\mu :=
\phi^+_\mu(\rho)^{4/n-2} \cdot \tilde g^\ast$  and $\mu \in [0,\mu_H]$: \:
$scal_{(\phi^+_\mu) ^{4/n-2} \cdot  \tilde g^\ast} \cdot (\phi^+_\mu) ^{4/n-2} \ge \iota_H/2\rho^2$.}\\

\textbf{Proof} \quad This follows readily from the transformation law
\[  - 4 (n-1)/(n-2) \cdot \Delta \phi^+_\mu  + scal_{g_H} \cdot \phi^+_\mu  =
  scal_{(\phi^+_\mu )^{4/n-2} \cdot
 g_H} \cdot (\phi^+_\mu) ^{n+2/n-2} \]
 Outside $C_p \setminus ( B_2(0) \setminus B_1(0))$, we just have $\Delta \phi^+_\mu = 0$ and in $C_p \cap ( B_2(0) \setminus B_1(0))$ we have
{\small \begin{eqnarray*} | \Delta \phi^+_\mu |
&= & | \Delta \varphi_\mu(\rho) \cdot \chi(\rho) + 2 <\nabla \varphi_\mu(\rho), \nabla  \chi(\rho)> +  (\varphi_\mu(\rho)+1) \cdot \Delta \chi(\rho) - \Delta \chi(\rho) |\\
&=& | 0 + 2 \mu \cdot < \nabla  \frac{1}{\rho^{n-2}}, \nabla  \chi(\rho)> + \mu \cdot \frac{1}{\rho^{n-2}} \cdot \Delta \chi(\rho) | \le const. \cdot \mu
\end{eqnarray*}}
\qed

Now we want to check that for $C  \in \overline {\cal T}_H$ equipped with $
\phi^+_\mu(\rho)^{4/n-2} \cdot \tilde g^\ast$ an $n-1$-dimensional area minimizer homologous to $\p B_1(0)$ stays outside a ball (whose radius depends on $\mu$) around $0$.\\

We will use this simple case as a sample how to handle the general problem that the ambient space is now singular (where no
general theory for area minimizers is available). \\
The method is a soft version of the whole approach: we can also bypass this problem by a doubling and choose a conical neighborhood $W_\gamma \subset V_\gamma$,
$1\gg \gamma$, with respect to $g$, such that $Vol_{n-1}\p (W_\gamma \cap \p B_1(0)) \le
\gamma$.\\
This can be done using the fact that the Hausdorff-dimension of $\sigma \cap \p B_1(0)$ is less than $n-3$: take a covering by balls $B_{\varsigma_i}(p_i)
\subset V_{\gamma/2}$ as in the definition of the Hausdorff-measure. For any $\zeta > 0$ we can find such a covering with $\sum_i \varsigma_i^{n-2} < \zeta$.
Thus for sufficiently small $\zeta$ one could define $\p W_\gamma$ (and thus $W_\gamma$) as a smoothed version of
$ \p \bigcup_i B_{\varsigma_i}(p_i) \subset \bigcup_i \p  B_{\varsigma_i}(p_i)$.\\

Now we delete $W_\gamma$ from $C$ and take a second mirror copy of $C \setminus W_\gamma$ and glue these two copies along there isometric boundary $\p
W_\gamma$. Except for the origin this gives a smooth manifold $D_\gamma$ with a $C^0$-metric. We ($C^\infty$-)smooth this metric symmetrically on $W_{\gamma +
\gamma^5} \setminus W_\gamma$ (and its mirrored copy) keeping the metric $C^0$-close to the original one. Consider the cone equipped with the metric
 $g_\mu$ running through such a doubling process also gives smooth metrics $C^0$-close to the respective doubled cone.
We call a suitably smoothed metrics $g^D_{\mu, \gamma}$ .\\

{\bf Lemma (6.5)}{\itshape \quad There is  some $\Theta_\mu \thicksim \mu^{1/n-2} >0$ such that for any cone $C$ there is a $\gamma_C >0$  such that for $0
<\gamma < \gamma_C $ the support of any area minimizing hypersurface $N^{n-1}$ in $(D_\gamma, g^D_{\mu, \gamma})$ homologous to the doubled unit sphere is still
outside the doubled sphere of radius $\Theta_\mu$. Moreover, for small radii the distance spheres will have $trA_{S_\rho(0)}(g)(v,w) \approx (n-1)\cdot \rho$ .} \\

 \noindent (Our sign convention is that $\p B_1(0) \subset \overline{B_1(0)}$ has \emph{positive} mean curvature, while $\p B_1(0) \subset \R^n \setminus B_1(0)$
 has \emph{negative} mean curvature)\\

Since the area minimizers can be assumed to be mirror-symmetric the volume-indifferent choice of the $W_\gamma$ allows us to paraphrase this (as a substitute
for the missing
solid theory) for the original cone:\\

{\itshape For any $C  \in \overline {\cal T}_H$ there is some $\Theta_\mu >0$ such that the support of any area minimizing hypersurface
$F^{n-1}$ in $(C, g_\mu)$ homologous to $\p B_1(0)$ is still outside $\p B_{\Theta_\mu}(0)$.} \\

In those cases (which includes prototypes of local barriers as considered below) where the boundaries of sufficiently narrow neighborhoods of the set ${\cal D}=
\overline{G} \setminus G$ in (6.2) form a foliation by isotopic submanifolds we can argue by showing that these submanifolds have \emph{positive mean
curvature}: the growth of the area elements of these submanifolds \emph{decreases} along the normals directed towards ${\cal D}$. (cf. [K] 1.9 which still
matches in our case without
Ricci curvature bounds)  \\
In particular it is only when we assemble these building blocks to global objects that we have to follow the longer route of using the individual barrier
properties to get an area minimizers with
obstacles in order to obtain such submanifolds with positive mean curvature.\\

The second fundamental form $A_L(g)$ of a submanifold $L$ with respect to some metric $g$ transforms under conformal deformations $g \ra u^{4/n-2} \cdot g$
according to following the formula $(AC)$
\[  A_L(u^{4/n-2} \cdot g)(v,w) = A_L(g)(v,w) - \frac{2}{n-2}\cdot {\cal{N}} (\nabla u / u) \cdot g(v,w)\]
where $ {\cal{N}} (\nabla u / u)$ is the normal component of $\nabla u / u$ with respect to $L$.\\

\textbf{Proof of (6.5)} \quad For the distance spheres $S_\rho(0)$ in $(D_\gamma, g^D_{\mu, \gamma})$ we have $\frac{2}{n-2}\cdot {\cal{N}} (\nabla \phi^+_\mu /
\phi^+_\mu) = - 2 \cdot \frac{\mu/\rho^{n-1}}{\mu/\rho^{n-2} + 1} \approx -2/\rho $ for $\rho \ll 1$  and for $1/\sqrt{\rho^{n-2}}$ (instead of $1/\rho^{n-2}$)
we get $\frac{2}{n-2}\cdot {\cal{N}} (\nabla \frac{1}{\sqrt{\rho^{n-2}}} / \frac{1}{\sqrt{\rho^{n-2}}}) = - \frac{1/\rho^{n/2}}{1/\rho^{{n-2}/2}} = -1/\rho $

Now consider $\left(1/\sqrt{\rho^{n-2}}\right)^{4/n-2} \cdot \tilde g^\ast = 1/\rho^2 \cdot \tilde g^\ast$ which is just a cylinder $\R \times S$. But in this
geometry $A_S \equiv 0$. Thus from $(AC)$: $trA_{S_\rho(0)}(g^\ast)(v,w) = -(n-1)/\rho$ with respect to $g^\ast$.\\

Applying $\varphi_\mu(\rho)$ and taking traces we therefore have for small $\rho$
\[ \phi^+_\mu(\rho)^{2/n-2} \cdot  tr A_{S_\rho(0)}(\phi^+_\mu(\rho)^{4/n-2} \cdot g)(v,w) \approx - (n-1)/\rho  + 2 \cdot (n-1)/\rho =   (n-1)/\rho\]
On the other hand $\phi^+_\mu = 1$ outside a ball and then the distance spheres are just the spheres within the cone and that means $tr A_{S_\rho(0)} = -
(n-1)/\rho$.\\
The latter makes sure that there is an area minimizer .\\
The way how the geometry changes (resp. reproduces) under scalings show the relation $\Theta_\mu \thicksim \mu^{1/n-2}$.\\

In order to get the desired result on the original cone $C$ we note that (by construction) for $\gamma \ra 0$ the volume of the set $W_\gamma$ will shrink to
zero and thus a potential area minimizer defined as a limit of Plateau problem (with free boundary data on the inner shrinking boundary) in the \emph{cone over}
$W_\gamma \setminus \sigma$ with boundary equal to $\p W_\gamma$ will eventually \emph{not} intersect $\frac{1}{2} \dot \p W_\gamma$ since otherwise we could
find after scaling a subsequence of these area minimizers converging to an area minimizer in a product with $\R$ (as a local blow up of a cone) whose support is
\emph{not} completely contained in the factor (perpendicular to $\R$) although the orthogonal projection is area contracting. \qed

\bigskip

Now we inductively describe how to define the \emph{local barriers} on a minimal cone equipped with the conformally deformed $scal > 0$-metric $\tilde g^\ast$ of (4.7). The starting point where the cone has only a singular tip has just been done.\\

Presuming we are done (with local and global barriers) for any case in any lower dimension we start again with carrying out an elementary deformation around the
tip $0 \in C$ such that the area minimizer stay
outside a ball of radius $4\rho$ (this forms the inner part of the local barrier).\\
Next we define the \emph{outer portions} of the local barriers. These parts are built inductively where the induction step is taking a lower dimensional global
barrier and consider a discretized product (with the usual scaling/tapering in cone direction) which covers a neighborhood of $\sigma \cap B_R(0) \setminus
B_\rho(0)$. For
given radii we can refine  the lower dimensional barrier such that its support is arbitrarily close to $\sigma$.\\

In the first step of isolated singular points comprising $\sigma$ we have $\sigma^* = \sigma$. The definition of $\sigma^*$ is also by induction and thus
becomes clear from the induction step \emph{plus} the construction of global barriers in the next section. Thus we may assume we have barriers
for  $\p B_1(0) \cap \sigma^* \subset \p B_1(0)$\\
\[\sum_{p \in P^{n-1}}\frac{ \lambda_p}{d_{\p B_1(0) \cap C}(x,p)^{n-3}} +1  \mbox { for some finite subset } P^{n-1} \mbox{ of } \p B_1(0) \cap \sigma^* \]
 and coefficients $ \lambda_p > 0$ such that the \emph{weights}
$\lambda_p$ and
the \emph{total weight} $\Theta \equiv \sum_{p \in P^{n-1}} \lambda_p$ are as small as ever needed (to ensure $scal > 0$ after truncation).\\

For the present we can suppress the truncation and handle the cut-off effect only when we estimate the negative effect on $scal$
below keeping in mind that they are of course needed to complete the definition when we install these deformations on $H$.\\

In order to define the \emph{outer part of the local barrier} we notice that the Laplacian for this product type space has to be calculated based on
$c(\omega)^{4/n-2} \cdot g_{\R} +
r^2 \cdot g_{\p B_1(0) \cap C}$ not on $ g_{\R} + r^2 \cdot g_{\p B_1(0) \cap C}$. Therefore (to keep $\Delta = 0$) we shift the exponent correspondingly, namely as follows:\\
$c(\omega)$ has the development $c(\omega) = \bar c_p (\bar \omega) \cdot r^{\beta_p}$ in $p$ (i.e. on a given tangent cone of $\p B_1(0) \cap C$ in $p$) and
because $(r^{\beta_p})^{2/n-2}$ is integrable on $\R^+$  the distance to $0$ and any other point in a cone $C$ conformally deformed by $(r^{\beta_p})^{4/n-2}$
remains finite. The distance of $\p B_t(0) \subset C$ to $0$ in the transformed metric is $(1 + \beta_p \cdot \frac{2}{n-2})^{-1} \cdot t^{1 + \beta_p \cdot
\frac{2}{n-2}}$. That means that $r = (1 + \beta_p \cdot \frac{2}{n-2}) \cdot \rho^{1/(1 + \beta_p \cdot \frac{2}{n-2})} $ for the transformed metric and thus
we have $(1 + \beta_p \cdot \frac{2}{n-2})^{2\beta_p/n-2} \cdot \rho^{2 \beta_p/(n + \beta_p -2)} $ as the (non squared) conformal factor in a point of distance
$\rho$.

Thus the $n$-dim deformation on $(B_1(0) \setminus B_\varsigma(0)) \cap C$ - which should be thought of as a Green's function along a 1-dimensional subspace =
cone-direction (within this non-smooth geometry) - to define the outer part of the local barrier is by taking a tapered product (where we suppress writing the
uniformly bounded constant $(1 + \beta_p \cdot \frac{2}{n-2})^{2\beta_p/n-2}$
\[(TP) \;\;\;\;  \sum_{p \in P^{n-1}}\frac{ \lambda_p}{d_{\p B_1(0) \cap C} \left( x/d(x,0),p \right)^{n-3-2 \beta_p/(n + \beta_p -2)}} + 1
 \]

In order to derive estimates we can reduce this to the case of products of the lower dimensional barrier with some compact interval since we can chop the
tapered version into pieces $B_{r+1/m} \setminus B_r \cap V_{1/m}(\sigma)$ (where $V_{1/m}$ is a cone shaped neighborhood of $\sigma$ of radius $1/m$ for $\p
B_1(0) \cap V_{1/m}$) for large $m \gg 1$ and scaled by $m^2$ and observe that (while these pieces blow up) the metric in any given point converges (from a
warped to a true) product metric in the sense that the corresponding of off-diagonal entries in the metric shrink to zero relative to diagonal entries and this
in
uniform way since we took already care of $\sigma \cap  B_{4\rho(0)}$ leaving a compact interval of radii in $\R^{>0}$. \\

Thus, for the calculus, we can focus on the following situation: $S \times \R := (\p B_1(0) \cap C) \times \R$ (with the warped product metric $g_S +
c(\omega)^{4/n-2} \cdot g_{\R}$ which can of course also be written $c(\omega)^{4/n-2} \cdot (g_{\R} + g_{\p B_1(0) \cap C})$ where the latter this time meant
the original metric on the tangent cone) and assume we already found a global barrier of truncated \emph{n-1 dim} Green's functions for an enhanced singular set
(to be defined in the proof) $S \cap \sigma^* \subset S$ of $S \cap \sigma \subset S $ in $(S, c(\omega)^{4/n-3} \cdot g_C|_{S})$, where we use that the Perron
property also applies to $c(\omega)$ for the dimensionally shifted eigenvalue equation $(CW)$ for  $S \cap \sigma^* \subset S$ and a finite $P^{n-1} $
 and coefficients $ \lambda_p > 0$ as above
\[ (P) \;\;\;\;  \sum_{p \in P^{n-1}}\frac{ \lambda_p}{d_{S \times \R}(x,\R \cdot p)^{n-3-2 \beta_p/(n + \beta_p -2)}} +1  \mbox {\;\;  for any } x \in S \times \R \]

Note that this product $S \times \R$ is an \emph{auxiliary object} useful for some  transparent computations which however is \emph{not} the tangent cone of any
point in $\sigma$. But this property is not used when the estimates for barriers close to $\p B_1(0) \cap \sigma^* \times \R$
are transformed into those for the (compactly supported) tapered case.\\

The computation in the proof of (6.3) gives again $\Delta = 0$ for $(P)$, since (on a tangent cone) we now  have $\sqrt{det(g_{ij},\{i \ge 2, j \ge 2\})} =
\sqrt{det(g_{ij})} = \rho^{n-2} \cdot (1 + \beta_p \cdot \frac{2}{n-2})^{2\beta_p/n-2} \cdot \rho^{2 \beta_p/(n + \beta_p -2)} $. \\

Our cones equipped with $\tilde g^\ast$ have $scal
>0$ in a scaling invariant fashion therefore we can use the blowing up view above to notice that on arbitrarily small neighborhoods of $\sigma$ the case $(TP)$ can
be handled in practice as the product from local approximations by $(P)$ (i.e. under scaling $\Delta \ra 0$ relative to the fixed $scal >0$).   \\

We continue with a look at the scalar curvature effect of a truncation of  $(P)$. The barrier property is handled in the next section.\\

We may assume by \emph{induction} that
\begin{enumerate} \item  a small neighborhood $V$ of $\p B_1(0) \cap C \cap \sigma$ in $\p B_1(0) \cap C$ can be chosen such that (applying (4.8))
$scal(\tilde g^\ast_{\p B_1(0) \cap C})> 0$ on $ V$
 \item the enhanced singular set $\sigma^\ast$ is also in $V$ i.e. $\p B_1(0) \cap \sigma^\ast \subset V$,
 \item there is a smaller neighborhood $W$ of $\sigma^\ast$ such that
\[\chi_{W} \cdot \sum_{p \in P^{n-1}} \frac{ \lambda_p}{d_{\p B_1(0) \cap C}(x,p)^{n-3}} +1  \mbox { for some finite subset } P^{n-1} \mbox{ of } \p B_1(0) \cap \sigma^* \]
where  $\chi_{W} \in  C^\infty(\p B_1(0) \cap C,[0,1])$  with  $\chi_{W} \equiv 1$ in $W$ and $\chi \equiv 0$ outside a small neighborhood of $W$ in $V$ is
\emph{barrier} for $\p B_1(0) \cap \sigma^\ast$,\\ Actually, the inductive construction uses a slightly more technical modification: we approximate $\sum_{p \in
P^{n-1}} \frac{ \lambda_p}{d_{\p B_1(0) \cap C}(x,p)^{n-3}}$ by a Riemann-Stieltjes sum of elementary barriers (see below)\item $scal((\chi_{W} \cdot \sum_{p
\in P^{n-1}} \lambda_p/d_{\p B_1(0) \cap C}(x,p)^{n-3} +1)^{4/n-3} \cdot
g^\ast|_V) > 0$ on $V$.\\
\end{enumerate}

{\bf Lemma (6.6)}{\itshape \quad Then we may also assume for (P) on $V \times \R$ (extending $\chi_{W}$ trivially in $\R$-direction) {\small \[scal(\chi_{W}
\cdot \sum_{p \in P^{n-1}} \lambda_p/d_{S \times \R}(x,\R \cdot p)^{n-3-2 \beta_p/(n + \beta_p -2)} +1)^{4/n-2} \cdot  (g_S + c(\omega)^{4/n-2} \cdot g_{\R}) >
0 \]}.}

\textbf{Proof} \quad Tracing back the definitions we see that both metrics on $\p B_1(0)$ and now on $\p B_1(0) \times \R$ had been deformed by the same factor
$c(\omega)^{4/n-2}$. Thus we can compare the two scalar curvature as in (6.4) via the transformation law (TL) and the fact that the interesting terms are \\
Now recall that under the transition $g \ra e^{2f} \cdot g$ on an \emph{n}-manifold we have $\Delta_gF \ra \Delta_{e^{2f} \cdot g}F = e^{-2f} \cdot (\Delta_gF +
(n-2)< \nabla f,\nabla F>)$.\\

Abbreviating $\sum_{p \in P^{n-1}} \lambda_p/d_{S \times \R}(x,\R \cdot p)^{n-3-2 \beta_p/(n + \beta_p -2)} +1$ by $\Psi$ we obtain \[| \Delta (\Psi \cdot
\chi_{W})| = | \Delta \Psi \cdot \chi_{W} + 2 <\nabla \Psi, \nabla \chi_{W}> + \Psi \cdot \Delta \chi_{W} | \] \[= \\ | 2 <\nabla \Psi, \nabla \chi_{W}> + \Psi
\cdot \Delta \chi_{W} | \le 2 |\nabla \Psi|\cdot |\nabla \chi_{W}| +  |\Psi \cdot \Delta \chi_{W} |\] and by comparing (TL) for $\chi_{W} \cdot \sum_{p \in
P^{n-1}} \frac{ \lambda_p}{d_{\p B_1(0) \cap C}(x,p)^{n-3}} +1 $ with the corresponding terms for $\chi_{W} \cdot \sum_{p \in P^{n-1}} \lambda_p/d_{S \times
\R}(x,\R \cdot p)^{n-3-2 \beta_p/(n + \beta_p -2)} +1$:
\begin{enumerate}
\item $| \nabla \chi_{W}|$ does not change under the dimensional ascent to (P). \item $|\Psi \cdot \Delta \chi_{W} |$ changes by $\Psi \cdot < \nabla c(\omega),
\nabla \chi_{W}>$
\end{enumerate}

We can already take care of these differences in the construction of $(\chi_{W} \cdot \sum_{p \in P^{n-1}} \lambda_p/d_{\p B_1(0) \cap C}(x,p)^{n-3} +1)^{4/n-3}
\cdot g^\ast|_V$: note the (only potentially problematic) cut-off effects are localized on a compact domain outside $\sigma$ which can be chosen inductively
(along with the the $P$ and the coefficients $\lambda_p$) in such a way that we can the coefficients for $| \nabla \chi_{W}|$ and $< \nabla c(\omega), \nabla
\chi_{W}>$ are that small that the change to $\Psi, \nabla \Psi$ is still small enough to be
 dominated be $scal > 0$. One iterates this reverse induction until one reaches single point singularities where (6.4)
applies. \qed

Using Riemann-Stieltjes sums one can see that $(P)$ (and via approximations also $(TP)$) can be considered as a compactly uniform limit for $l \ra \infty$ on
$(B_1(0) \setminus B_\varsigma(0)) \setminus \sigma^*$ of a sequence of sums of finitely many \emph{n-dimensional} Green's functions locally defined on well
$C^k$-approximated parts of balls  around points $(p, t )$ by some tangent cones $C_{(p, t )}$: $\varphi_{(p, t )}$.  As already mentioned we will only need to
care for a small shell around the singular sets (though neither the fact that these functions are not defined globally nor near a worse approximated part of $H$
matters).\\

Thus up to some normalizing constant $(P)$ is a limit for $l \ra \infty$ of
\[ \sum_{t\in 1/l \cdot \Z } \sum_{p \in P^{n-1}}
\frac{\lambda_p}{l} \cdot \varphi_{(p, t )}(x)) +1 \] Also we get a way to truncate $\sum_{p \in P^{n-1}} \frac{ \lambda_p}{d_{\p B_1(0) \cap C}(x,p)^{n-3}}$ by
a sum of truncated functions (the individual negative cut-off effects for $scal$ are still captured from (6.4) since sum remains the (at most the same) since
the coefficients enter linearly in these estimates. That is the outer portion of a local barrier (modulo the approximation of (TP) by (P) is
\[ \sum_{t\in 1/l \cdot \Z } \sum_{p \in P^{n-1}} \frac{\lambda_p}{l} \cdot \phi_{(p, t )}(x) +1 \]

 In this paper this finite sum (for large $l$) is merely used for its expositorily advantages,
there are two purposes:
\begin{enumerate} \item the inductive definition based on barriers in lower dimensions always consisting of finitely balls/Greeen's functions as above,
\item in order to extend this to domains outside $(B_1(0) \setminus B_\varsigma(0)) \cap C$ just choosing the coefficients $=0$, in place of an artifical
cut-off construction we had to introduce for $ \lambda_p/d_{\p B_1(0) \cap C} \left( x/d(x,0),p \right)^{n-3-2 \beta_p/(n + \beta_p -2)}$.
\end{enumerate}
Thus we do not present this computation, instead we notice a simpler argument to see why such a presentation exists: each of these functions and hence their sum
has $\Delta = 0$ and we get with some normalization along some boundary of  neighborhoods of the $p$ that these sums converges to $(P)$ along this
boundary and this leads to that compact interior convergence.\\

\vspace{1.7cm}
\begin{center}
\large{\textbf{7. Global barriers}}
\end{center}

Now we will we assemble the local and global barriers and verify that they have the desired deflection properties.\\

We start with the definition of an entire \emph{local barrier} on a tangent cone $C$, we first take:
\[\left(\frac{\lambda_0}{d(x,0)^{n-2}} + \sum_{p \in P^{n-1}}\frac{ \lambda_p}{d_{\p B_1(0) \cap C} \left( x/d(x,0),p \right)^{n-3-2 \beta_p/(n + \beta_p -2)}} + 1\right)^{4/n-2} \cdot \tilde g^\ast\]\\
For $\lambda_0 \ra 0$ we know that the distance of the  area minimizer homologic to $\p B_1(0)$ closest to $0$ is  $\approx  \kappa_0 \cdot
\lambda^{1/{n-2}}_0$.\\ For  small $\rho$ we  observe an overlap of the barriers for $\left(\lambda_0/d(x,0)^{n-2}   + 1\right)$ and $\left(\sum_{p \in P^{n-1}}
\lambda_p/d_{\p B_1(0) \cap C} \left( x/d(x,0),p \right)^{n-3-2 \beta_p/(n + \beta_p -2)} +1\right)$  leading by the previous and the argument below to a common
barrier if the \emph{weights} $\lambda_p$ and
the \emph{total weight} $\Theta \equiv \sum_{p \in P^{n-1}} \lambda_p$ are chosen small enough and $P$ is chosen suitably.\\

And again to complete the definition of the local barrier we will use the Riemann-Stieltjes sum presentation with zero coefficients for contributions outside
the distance balls in the original cone metric (used for the covering argument) that allows us to restrict the deformation to a ball in
a natural way.\\

This construction of local barriers in the \emph{cone} case can directly be transferred (by some push-forward) to balls around singular points in $H$  since the
approximation of $H$ is arbitrarily fine in $C^k$-topology on sets $(B_{R}(0) \setminus B_{\varrho}(0)) \setminus V_{\xi}(\sigma)$, cf. [L5] pp. 667 - 669 for
technical details.\\
The trick is that those parts not well approximated (in the respective ball) will be hidden behind the barrier.
Hence, locally, the cone barrier construction transplants to $H$.\\

After these preparations we eventually turn to the barrier property.\\

There are two cases we have to understand: the first one is the dimensional transition property of lower to higher dimensional barriers, that is, we will check
that (tapered) products of lower dimensional global barriers of the previous section (P) and (TP) lead again to barriers.\\ The second one below (which also
covers the composition of inner and outer part of a local barrier) is the property that the superposition of an upper bounded number
of such barriers for certain $\sigma^\ast$'s forms a barrier for a union of these $\sigma^\ast$'s. \\

Both operations lead to a perturbation of the barrier property and we overcome this using extreme  choices of barriers very
close to $\sigma^\ast$ and thus (after the deformation) far outside some kind of horizon although (for a fixed deformation) not arbitrarily close to infinity
since eventually the single point set up of the obstacle becomes visible again.\\

We start inductively with an isolated point singularity in a cone $C^{n-1}$. Then we observe that we get considering distance tubes instead of spheres for $C
\times \R$  from the computation in (6.4) that \[(1/d_{C \times \R}(x,\R \cdot p)^{n-3-2 \beta/(n + \beta -2)} +1 )^{4/n-2} \cdot c( \omega) ^{4/n-2} \cdot
(g_C+ g_{\R})\] also forms a barrier. Now we cut this off considering for $l \gg 1$
\[ \sum_{t\in 1/l \cdot \Z \cap [0,1]}
\frac{1}{l} \cdot \varphi(x,t) +1 \mbox{ for } (x,t) \in C \times \R \] As in the proof of (6.5) we observe that this forms an obstacle for $\sigma \times
[0,1]$ for $l \ra \infty$.\\
A  straightforward extension of this argument applies to the
(warped)
product of higher dimensional barriers with a copy of $\R$ respectively an interval.\\

 Notice that the outward bending effect towards $0$ decreases when we raise the dimension but this is compensated
from the bending effect for the additional $\R$ factor.
 (For the surgery type argument this corresponds to the fact that we cannot expect to reach the situation of totally geodesic boundaries just using conformal
deformations.)\\

We can now superpose 1-dimensional (and inductively higher
dimensional) local barriers on a minimal hypersurface to get our
\emph{global barrier} for its singular set. This is done with the
covering tools from sec 3. to find $c(n)$ families of local
barriers where barriers in a given family do not intersect and
superpose the local barriers, that is the sum of all truncated
Green's functions defined on one tangent cones on each $B_{
\varrho_{ij}}(p_{ij}) \in {\cal{F}}^a$.\\
Note that this will
inevitably be a barrier not just for the singular set but a 1- and
inductively  higher dimensional stratified approximation of the
singular set in the following sense:\\

\textbf{Definition (7.1)}{\itshape \quad  For an area minimizing
hypersurface $H$ with singular set $\Sigma$ we define inductively
an \emph{\textbf{enhanced singular}} set $\sigma^\ast$  as
\[\bigcup_{B_{ \varrho_{ij}}(p_{ij}) \in {\cal{F}}^a}B_{
\varrho_{ij}}(p_{ij}) \cap \sigma^\ast_{p_{ij}} \mbox{ with the
summands } B_{ \varrho_{ij}}(p_{ij}) \cap \sigma^\ast_{p_{ij}}
\subset C_{p_{ij}}\] where $\sigma^\ast_{p_{ij}}$ means a
\emph{\textbf{crystallized}} lower dimensional enhanced singular
sets for $\sigma_{p_{ij}} \cap \p B_1(0) \cap C_{p_{ij}}$.

We call a subset of $H$ (which one may imagine as $\Sigma \: \cup$
\emph{some highly curved regions in $H$}) a
\textbf{\emph{crystallization}} of $\sigma^\ast$ if it is near to
$\sigma^\ast$ in the sense of $\Sigma \prec \sigma$ (for each
involved tangent cone and ball) and the barrier built using
$\sigma^\ast$ can be
written as a Riemann-Stieltjes sum of single Green's functions centered in points of this crystallization.}\\

 For the induction step (where we suppose to have a fine approximation of cones on these balls) each $\sigma^\ast$ is crystallized and henceforth identified with
 its crystallization since we only
care for a shell around this part the balls are only formally
assigned to points in $H$ while the choice came from points in
$\sigma^\ast$. In other words this
very unsharp notion is enough since we use it by its functionality.  \\

 It is readily seen
that the under composing the global barrier (for suitably chosen
parameters) the scalar curvature remains $> 0$: since each of the
local barriers keeps it positive we can adjust the coefficients
along (6.4) to allow the $c(n)$-fold perturbation
without losing positivity.\\

\textbf{Proposition (7.2)}{\itshape \quad For suitably chosen
parameters this makes up a global barrier for $\Sigma$}.\\

{\bf Proof}\quad In order to understand the effect of a superposition of barriers the sharpness of approximation (in the sense of the relation $\Sigma \prec \sigma$) becomes important:\\
If we knew that the approximation was perfect (= the
crystallization is just the isometric identification) then the
superposition would quite obviously give a global barrier since
the poles are neither moved nor damped by a superposition of two
(and hence by a universally upper bounded number of)
obstacles and therefore they were surrounded by a merged obstacle.\\
This is also a good place to see why the intermediate creation of local barriers (instead of using elementary barriers directly to go for a global barrier) was
important: even in such a perfect situation elementary barriers definitively secure only a single point from being hit by a minimizer. Every other point could
be
secured just by a fortunate choice of other elementary barriers (and a local barrier is such a choice).\\

 Thus one uses the covering tool to provide a uniform approximation (modulo scaling) to come sufficiently close to this situation.  Although we used this already several times  rather silently this is a place
 to mention that refined approximations etc. means that we choose finer coverings ${\cal{F}}^a$ for $a \gg 1$ and the statements are valid with the reservation that
 we have chosen a sufficiently large $a$.   \\

As a technical aid we can now use minimal hypersurfaces with an
obstacle compressing around the potential barriers to localize
this problem. This is defined for a local barrier on the cone
under consideration as a tube around $\sigma^\ast$ such that the
boundary is in the well $C^k$-approximating domain. The boundary
of the union of the interior of these tubes (truncated by the size
of the respective ball) will be called $W^+$ and we take a
corresponding $W^-$ as the supposed not reached inner obstacle
(leaving the parameters open for adjustments).\\
Recall also that the minimizer in this homology class would be a point.\\

The proof is by induction over $k$, that is, we activate all local
obstacles in the \emph{k-th} family ${\cal{F}}_k$
\emph{additionally} to those in the families ${\cal{F}}_l, \; l <
k$ and we
will check their common deflecting properties.\\
 Since the obstacles within one family stay separated the case $k=1$ is just the local barrier property. \\
The induction assumption is that for the local barriers in $\bigcup_{m \le k} {\cal{F}}_k$, $k < c(n)$ and $(W_k^+, W_k^-)$ form a barrier for the union of the
transplants to $H$ of the respectively scaled (by $\varrho_{ij}$) versions in the chosen tangent cone in $p_{ij}$.\\

 In turn this allows us to find a localization
effect of $W^-$ and $W^+$ as obstacles compressing an area
minimizer $N^{n-1}$ homologous to $\p W^+$ in the difference
region $W^+ \setminus W^-$: the mutual influence on the shape of
the area minimizer next to different members within $c(n)$
disjoint families ${\cal{F}}_l \subset {\cal{F}}^a$ becomes
arbitrarily small.\\

Hence we can localize the induction argument starting with the assumption that the hypothesis was proved for any combination of some $k < c(n)$ local barriers
of $k+1$ given local barriers (belonging to different families ${\cal{F}}_l$) and we are done if we can prove that the combination of all these  $k+1$ local
barriers also form a
barrier for the respectively larger union of enhanced pieces of $\Sigma$.\\

Moreover, we may restrict to the case where we have only two local barriers since iterating this argument a uniformly upper bounded number of times can be
handled by choosing the ratio of sizes correspondingly larger which allows to improve
the estimates as needed (cf. [L5],pp 671-672 for the way how to carry out such an argument formally):\\

Using these preparations we can finally reduce the discussion to the case of two local obstacles and actually we can take two outer parts of different barriers
(the union of outer and inner part of a local barrier can be handled considering the inner part as the outer one of another barrier).\\
Moreover the overlap of the two barriers takes place on the intersection of the two carrier balls closely approximated by cones. Since the single barriers are
in the smooth fine approximation parts we can argue using the compressing obstacles that we can carry this question out on the cones. If we get a new barrier
and it is in the well-approximated part this can be taken over for $H$.\\

But here we just notice (as already stated) that in the case of a
perfect approximation we find such a inner unreachable obstacle
$W^-$ and since we are free to refine the cone approximation to
any extend this obstacle will eventually be in the well $C^k$-
approximated portion and thus remains a barrier for the area
minimizer

 \qed

 \vspace{1.7cm}
\begin{center}
\large{\textbf{8. Obstacles leading to totally geodesic borders}}
\end{center}

\bigskip

Since $(W^-,W^+)$ forms a barrier for $\Sigma$ we know that the
area minimizer $N^{n}$ in this region no longer reaches $\p W^-$
(and hence $\Sigma$) but it may touch some parts of $\p W^+$. We
make a virtue of necessity and modify $\p W^+$: cover $\p W^+$ by
finitely many small open balls $B(x_i)$ for some $x_i \in \p W^+$
such that the mean curvature of every $\p B(x_i)$ is
\emph{strictly positive} and such that $W^+_\flat = W^+ \setminus
\bigcup_i B(x_i) $ can be used
instead of $W^+$, that is, $W^+_\flat \supset W^-$ and $(W^-,W^+_\flat)$ forms a barrier for $\Sigma$.\\
If we take an area minimizer with these obstacles we cannot say
much about its regularity, but we push this one step further and
deliberately use such contact balls to gain global regularity
while loosing the vanishing mean curvature: after placing
additional balls suitably we conclude that the resulting
hypersurface
is $C^1$.\\

{\bf Proposition (8.1)}{\itshape \quad Let $V^{n}$ be the unique area minimizer with obstacles and within some homology class of some compact orientable
$(M^{n+1},g_M)$.\\ Assume that all obstacles are in contact on one-side of $V^{n}$, then we can place additional obstacles on the same side and find a new
$C^1$-smooth area minimizer ${\cal{V}}^{n}$ arbitrarily near to $V^{n}$ in Hausdorff-topology within the same
homology class with both classes of objects as obstacles.}\\

The proof of this result is technically unrelated to the present paper and appears in [L6]. \\

This $C^1$-smooth area minimizer ${\cal{V}}^{n}$ (with obstacles)
close to $V^{n}$ has the feature that in those places of
coincidence with the boundaries of the balls it has positive mean
curvature and otherwise (as a free minimal surface) zero mean
curvature. Next we can smooth ${\cal{V}}^{n}$
to get a smooth ${\cal{W}}^{n}$ close to $V^{n}$ with \emph{positive mean curvature} (cf. [L6]).\\

Finally we deform an arbitrarily small neighborhood of
${\cal{W}}^{n}$ making it \emph{totally geodesic} and such that
the
scalar curvature on $H$ \emph{increases }in those parts where we change the metric. \\
This is based on an \emph{h-principle} type consideration of the
curvature expression and uses the linearity in second derivatives
while it is homogenous quadratic in first order terms.
Geometrically, we compress geodesics which are leaving the
boundary perpendicularly.\\

\textbf{Proposition (8.2)} \emph{For $U$ be a neighborhood of ${\cal{W}}^{n} \subset H^n .$ Then, there is a metric $g_{U}$ on $H$ with the following
properties:}

\begin{enumerate}
\item $scal (g_{U}) \ge scal (g)$ \item ${\cal{W}}^{n}$ \emph{is
totally geodesic with respect to} $g_{U}$  \item $g_{U} \equiv g$
\emph{on the regular component of} $H \setminus U$
\end{enumerate}
\medskip \emph{Note that the proof of this result also concludes
the \textbf{proof of our Theorem} in the introduction.}\\
\medskip

In order to construct $g_{U}$ note that we mainly have to care
about a one sided neighborhood $V$ of ${\cal{W}}^{n}$  within $H$
since we will delete the part that is on the $\Sigma$-side. The
metric $g$ can be written as $g_{\R} + g_t$ when identified
isometrically with $\left( [0, \sigma) \times  {\cal{W}}^{n},
g_{\R} + g_t \right),$ for $\sigma$ small enough and $g_t$ defined
as $g|_{TN_t}$ on $N_t =\{ x \in H| dist (x, {\cal{W}}^{n}) = t
\}$, since the geodesics orthogonally emanating from
${\cal{W}}^{n}$ hit each $N_t$ orthogonally.
\\

Now the idea is to substitute the metric $g_{\R} + g_t$ on $[0,
\sigma) \times {\cal{W}}^{n}$ for $g_{\R} + g_{h (t)}$ where $h
(t)$ is a function with strongly positive
second derivative (when compared with its first derivative): \\

\textbf{Lemma (8.3)} \quad \emph{For any given for $k > 0$ there is a smooth function $h
> 0$ on $\R$ with $h(t) = - t$ on $\R^{\le -\delta}, h(t) = t$ on $\R^{\ge \delta}, h (-
s) = h (s)$, for some $\delta \in (0, \sigma/2),\: |h' (r)| \le 1, \: h^{''} \le 0$ such that $|h^{''}(r)| \ge k \cdot |h' (r) + 1 |$ on $\R^{\le 0}$ and
$|h^{''}(r)| \ge k
\cdot |h' (r) - 1 |$ on $\R^{\ge 0}$.} \\

This is an elementary construction combining functions of the type $f (x) = s \cdot \exp (-d/t)$ with cut-off functions. Now we claim\\

\textbf{Proposition (8.4)} \quad \emph{For sufficiently large $k$ we have}

\begin{enumerate}
\item $scal (g_{\R} + g_{h(t)}) \ge scal (g_{\R} + g_t)$, \item
\emph{and with respect to this metric} ${\cal{W}}^{n}$
\emph{becomes (obviously) totally geodesic.}
\end{enumerate}

\textbf{Proof} \quad For $g$ written in local coordinates $x_i : \sum_i g_{ij} dx_i dx_j$ we have the standard formulas { \small \[ scal (g) = \sum_{i,j,k}
g^{ij} \left( \frac{\partial \Gamma^k_{ij}}{\partial x_k} - \frac{\partial \Gamma^k_{ik}}{\partial x_j} + \sum_{l} \Gamma^l_{ij} \cdot \Gamma^k_{kl} -
\Gamma^l_{ik} \cdot \Gamma^k_{jl} \right) , \]
\[ \Gamma^\gamma_{\alpha \beta} = \frac{1}{2} \sum_p g^{\gamma p} \left( \frac{\partial
g_{\alpha p}}{\partial x_\beta} + \frac{\partial g_{\beta p}}{\partial x_\alpha} - \frac{\partial g_ {\alpha \beta}}{\partial x_p} \right), \]
\[ \frac{\partial \Gamma^k_{ii}}{\partial x_k} = Q (g, \frac{\partial g}{\partial x}) + \frac{1}{2} \sum_p g^{k p} \cdot
\left( \frac{\partial^2 g_{i p}}{\partial x_i \partial x_k} + \frac{\partial^2 g_{i p}}{\partial x_i
\partial x_k} - \frac{\partial^2 g_{i i}}{\partial x_p \partial x_k} \right) ,\]
\[ \frac{\partial \Gamma^k_{ik}}{\partial x_i} = Q (g,\frac{\partial g}{\partial x}) + \frac{1}{2} \sum_p g^{k p} \cdot
\left( \frac{\partial^2 g_{i p}}{\partial x_k
\partial x_i} + \frac{\partial^2 g_{k p}}{\partial x_i \partial x_i} - \frac{\partial^2 g_{k i}}{\partial x_p \partial x_i} \right)
\]}

where the \emph{Q}'s are quadratic terms in 1st order, rational in zeroth order derivatives of the metric.\\

We choose \emph{Fermi-coordinates} $x_1, \cdots x_n$ along
${\cal{W}}^{n}$ such that $x_1$ corresponds to the parameterized
normal geodesics emanating ${\cal{W}}^{n} $. Also we choose the
coordinates such that they are \emph{geodesic coordinates} on
${\cal{W}}^{n} $ around a base point $p \in {\cal{W}}^{n}$.

In these coordinates the growth $\frac{\p g_{mm}}{\partial x_1}$
is directly described by Jacobi fields. Observe that $x_1$ is the
first coordinate in the decomposition $g_{\R} + g_t$ on $[0,
\sigma) \times {\cal{W}}^{n}$  and $x_2, ... x_n$ are the
coordinates along $N_t$. In particular \[(A) \;\;\;\;\;\;\;\;
g_{1i} \equiv \delta_{1i},\;\; \frac{\p g_{1i}}{\partial x_p}
\equiv 0, \;\; \frac{\partial^2 g_{1i}}{\partial x_p \partial x_m}
\equiv 0 \mbox{ along } {\cal{W}}^{n}
\]  \[ \mbox{  and }  g_{ij} = \delta_{ij}, \frac{\p
g_{ij}}{\partial x_k} = 0 \mbox{ and therefore } g^{ij} =
\delta_{ij} \mbox{ in } p \]

Although these equations hold only on ${\cal{W}}^{n} $ resp. just
in $p$, this allows us to use estimates for these quantities in
small neighborhoods within $H$ and $K$.\\

Also, we note that precisely the terms $g^{ii} \cdot \frac{\partial \Gamma^1_{ii}}{\partial x_1}$ contain those contributions which will turn out to dominate
the curvature expression. Abbreviating the terms with mixed derivatives 2nd and quadratic 1st order derivatives by $ ( \cdot) $ we have { \footnotesize \[ (B)
\qquad g^{mm} \cdot \frac{\partial \Gamma^1_{mm}}{\partial x_1} = ( \cdot) - g^{mm} \cdot \frac{1}{2} \cdot \frac{\partial^2 g_{mm}}{\partial x^2_1} \cdot
g^{11}, \;\;\;\;\; g^{11} \cdot \frac{\partial \Gamma^m_{1m}}{\partial x_1} = ( \cdot ) + g^{11} \cdot \frac{1}{2} \cdot \frac{\partial^2 g_{mm}}{\partial
x^2_1} \cdot g^{m m} \]}

Note that the latter terms $g^{11} \cdot \partial \Gamma^m_{1m} /
{\partial x_1} $ appear with a \emph{minus} sign in \emph{scal }
and that for any $\gamma> 0$ we can find a $\sigma > 0$ such on
$[0, \sigma) \times {\cal{W}}^{n} $ :
\[ (C \gamma) \qquad 1- \gamma \le g^{ii} \le 1 + \gamma . \]

The main point is that for the deformation $g_{\R} + g_t \rightarrow g_{\R} + g_{h(t)}$ above the terms $-\frac{1}{2} \cdot \frac{\partial^2 g_{mm}}{\partial
x^2_1} \cdot g^{mm} \cdot g^{11}, m \ge 2$ (for $m = 1$ it vanishes since $g_{11} \equiv 1$) literally dominate the deviation of all the
other second, first and zeroth order terms appearing in \emph{scal}:\\

For the Jacobi field $Y (t)$ with $Y (0) =
\frac{\partial}{\partial x_m} \in T_p {\cal{W}}^{n} $ along the
normal geodesic $\gamma_{\nu}$ with $\gamma_{\nu} (0) = p$ and $Y'
(0) = A Y (0)$ we get :

\[(D) \qquad 2 \cdot A (\frac{\partial}{\partial x_m},
\frac{\partial}{\partial x_m}) =2 \cdot | Y (0) |^2 \cdot A
(\frac{\partial}{\partial x_m}, \frac{\partial}{\partial x_m}) /
g_{mm} (0)\] \[ = 2 \cdot \left< Y (0), Y' (0) \right> = \left( |
Y (t) |^2 \right)' (0) = \frac{\partial g_{mm}}{\partial x_1}|_{t
= 0}
\]


and this eventually leads to estimates for \emph{scal}:\\

\textbf{Lemma (8.5)} \quad \emph{For any $\ve > 0$ the reparameterization there is a $k$ such that for the associated $h(t)$:} {\small \[ | ( scal(g_{\R} + g_t)
- ( \sum_m -\frac{1}{2} \cdot \frac{\partial^2 (g_{\R} + g_t)_{mm}}{\partial x^2_1} \cdot (g_{\R} + g_t)^{mm} \cdot (g_{\R} + g_t)^{11} )  ) - \] \[\ (
scal(g_{\R} + g_{h(t)}) - ( \sum_m -\frac{1}{2} \cdot \frac{\partial^2 (g_{\R} + g_{h(t)})_{mm}}{\partial x^2_1} \cdot (g_{\R} + g_{h(t)})^{mm} \cdot (g_{\R} +
g_{h(t)})^{11} ) ) | \; \; \; \le
\] {  \[ \ve \cdot |  \sum_m \frac{1}{2} \cdot \frac{\partial^2 (g_{\R} + g_t)_{mm}}{\partial x^2_1} \cdot (g_{\R} + g_t)^{mm} \cdot (g_{\R} + g_t)^{11}  - \]
\[  \sum_m \frac{1}{2} \cdot \frac{\partial^2 (g_{\R} + g_{h(t)})_{mm}}{\partial x^2_1} \cdot (g_{\R} + g_{h(t)})^{mm} \cdot (g_{\R} + g_{h(t)})^{11}  | \]}}

\textbf{Proof} \quad We will actually show that such an inequality
holds for all the individual terms which appear in \emph{scal}.
We start with those terms containing only zeroth and first order derivatives of the metric:\\

There are of the two types: \: $ g^{i_1 i_2} \cdot \frac{\partial g_{i_3 i_4}}{\partial x_{i_5}} \cdot
 \frac{\partial g_{i_6 i_7}}{\partial x_{i_8}}$ \: and \: $g^{i_1 i_2}
 \cdot g^{i_3 i_4} \cdot g^{i_5 i_6} \cdot \frac{\partial g_{i_7 i_8}} {\partial x_{i_9}}
 \cdot \frac{\partial g_{i_{10} i_{11}}}{\partial x_{i_{12}}}$ \\

Since $g^{ij} = P_{n} (g_{v w}) / \det (g_{\alpha \beta})$, where
$P_{n}$ is a homogenous polynomial of degree $n$ in entries of
$(g_{\alpha \beta})$, we also note $\frac{\partial
g^{ij}}{\partial x_k} = \sum \frac{\partial g_{ab}}{\partial x_k}
\cdot P_{2n} (g_{v w}) \cdot \det (g_{\alpha \beta})^{-2}$ and get
that the two types of terms have the same elementary form:

\[(\ast) \qquad \frac{\partial g_{a b}}{\partial x_c} \cdot \frac{\partial g_{i j}}{\partial x_k} \cdot P_{3 n } (g_{v w}) \cdot \det
(g_{\alpha \beta})^{-3}
\]

Secondly, we also have terms containing second order derivatives
of the type \\
$g^{i_1 i_2} \cdot g^{i_3 i_4} \cdot \frac{\partial^2 g_{i_5 i_ 6}}{\partial x_{i_7}
\partial x_{i_8}}$, with $i_7$ or $i_8$ unequal $1$,
since we subtract precisely those terms with $i_7 = i_8 = 1$. \\
We observe that this kind of expressions is of the form \\
\[(\ast\ast) \qquad \frac{\partial^2 g_{ij}} {\partial x_k \partial x_m} \cdot P_{2 n } (g_{v w}) \cdot \det (g_{\alpha \beta})^{-2} \] Here
we will distinguish the two cases where $k = m = 1$ which plays the pivot role in the curvature
estimates and the case where $k \neq 1$ or $m \neq 1$ \\

Now we compare the various derivatives for $g_{\R} + g_{h (t)}$ with those of $g_{\R} + g_t:$ \\
First of all: $\left(g_{\R} + g_{h (t)}\right)_{a b} \Big{\vert}_{t = t_0} = \left(g_{\R} + g_s \right)_{a b} \Big{\vert}_{s = h (t_0)}$ and for $i, j \ge 2$
and $a, b \ge 1$ we obviously have
\begin{center}
$ \frac{\partial \left(g_{\R} + g_{h (t)} \right)_{a b}}{\partial x_i} \Big{\vert}_{t = t_0} = \frac{\partial (g_{\R} + g_s)_{a b}}{\partial x_i} \Big{\vert}_{s
= h (t_0)}, \frac{\partial^2 \left(g_{\R} + g_{h (t)}\right)_{a b}}{\partial x_i
\partial x_j} \Big{\vert}_{t = t_0} = \frac{\partial^2 \left(g_{\R} + g_s \right)_{a b}}{\partial x_i \partial x_j}
\Big{\vert}_{s = h (t_0)}$
\end{center}

\noindent The coordinates on $\left( [0, \delta) \times
{\cal{W}}^{n} , g_{\R} + g_{h (t)} \right)$ are defined to equal
those of $g_{h (t)}$ in
$\left([0, \delta) \times {\cal{W}}^{n} , g_{\R} + g_t \right)$. \\
This means if we want to understand $scal (g_{\R} + g_{h (t)})$ - better to say show $scal \ge 0$ - on the distance level surface $N_t$ we can use most of the
terms for $scal (g_{\R} + g_t)$ on $N_{h (t)}$ substituting only those terms $\frac{\partial g_{a b}}{\partial x_1}, \frac{\partial^2 g_{a b}}{\partial x_1
\partial x_c}, \mbox{ with } a, b \ge 2, c \ge 1$ for the corresponding ones in $scal (g_{\R} + g_{h (t)})$. These terms can be evaluated by elementary means:
$(F)$

\[ \frac{\partial}{\partial x_1} \left(g_{\R} + g_{h (t)}\right)_{ab} \Big{\vert}_{t = t_0} \;\; = \;\; h' (t) \cdot \frac{\partial}{\partial
x_1} \left(g_{\R} + g_s \right)_{a b} \Big{\vert}_{s = h (t_0)},\] \[ \frac{\partial}{\partial x_1} \left( \frac{\partial}{\partial x_c}\left(g_{\R} + g_{h (t)}
\right)_{a b}\right) \Big{\vert}_{t = t_0} \;\; = \;\; h'(t) \cdot \frac{\partial^2 (g_{\R} + g_s)_{a b}}{\partial x_1
\partial x_c} \Big{\vert}_{s = h (t_0)}, \ c \ge 2\]
{\footnotesize \[ \frac{\partial^2}{\partial x^2_1} \left(g_{\R} + g_{h (t)}\right)_{a b} \Big{\vert}_{t = t_0} \;\; = \;\; h'' (t) \cdot
\frac{\partial}{\partial x_1} \left(g_{\R} + g_s \right)_{a b} \Big{\vert}_{s = h (t_0)} + \left(h' (t)\right)^2 \cdot \frac{\partial^2}{\partial x^2_1}
\left(g_{\R} + g_s \right)_{a b} \Big{\vert}_{s = h (t_0)} \]}
\\
Thus we substitute those components of $scal \left(g_{\R} + g_t \right)$ containing one of these terms for their multiples by $h' (t)$, $h' (t)^2$ resp. $
h''(t)$ passing thereby to
{\footnotesize  $scal \left(g_{\R} + g_{h (t)}\right)$ \\
\[(\ast \ast \ast) \qquad scal \left(g_{\R} + g_{h (t)} \right)\Big{\vert}_{N_t} - scal \left(g_{\R} + g_t \right) \Big{\vert}_{N_{h (t)}} = \]
\[\sum\limits_{I_1} (h' (t)^2 - 1) \cdot \frac{\partial}{\partial x_1} \left(g_{\R} +
g_s \right)_{ab} \cdot \frac{\partial}{\partial x_1} \left(g_{\R}
+ g_s \right)_{ij} \cdot P_{3n} \left( (g_{\R} + g_s)_{vw} \right)
\cdot \det \left( (g_{\R} + g_s)_{\alpha \beta} \right)^{-3}
\Big{\vert}_{s = h (t)} + \]
\[\sum\limits_{I_2,k \geq 2} (h' (t) - 1) \cdot \frac{\partial}{\partial x_1} \left(g_{\R} +
g_s \right)_{ab} \cdot \frac{\partial}{\partial x_k} \left(g_{\R}
+ g_s \right)_{ij} \cdot P_{3n} \left( (g_{\R} + g_s)_{vw} \right)
\cdot \det \left( (g_{\R} + g_s)_{\alpha \beta} \right)^{-3}
\Big{\vert}_{s = h (t)}\]
\[+ \;\; \sum\limits_{I_3} (h' (t)^2 - 1) \cdot
\frac{\partial^2}{\partial x_1
\partial x_1} \left(g_{\R} + g_s \right)_{ab} \cdot P_{2n}
\left( (g_{\R} + g_s)_{vw} \right) \cdot \det \left( (g_{\R} + g_s)_{\alpha \beta} \right)^{-2} \Big{\vert}_{s = h (t)}\]
\[+ \;\; \sum\limits_{I_4,c \geq 2} (h' (t) - 1) \cdot
\frac{\partial^2}{\partial x_1
\partial x_c} \left(g_{\R} + g_s \right)_{ab} \cdot P_{2n}
\left( (g_{\R} + g_s)_{vw} \right) \cdot \det \left( (g_{\R} +
g_s)_{\alpha \beta} \right)^{-2} \Big{\vert}_{s = h (t)}\] {  \[+
\;\; \sum\limits_{I_5} h'' (t) \cdot \frac{\partial}{\partial x_1}
\left(g_{\R} + g_s \right)_{ab} \cdot P_{2n} \left( (g_{\R} +
g_s)_{vw}) \right) \cdot \det \left( (g_{\R} + g_s)_{\alpha \beta}
\right)^{-2} \Big{\vert}_{s = h (t)}\]}} where the restriction of
$scal$ to $N_t$ resp. $N_{h (t)}$ corresponds to the
identification of the corresponding points on these level
surfaces, and where $I_1, I_2, I_3, I_4, I_5$ are the respective
index sets covering those term which show up in $scal$.
\bigskip

In order to start playing with these entities note that, having
chosen Fermi coordinates about some point $p \in {\cal{W}}^{n} $,
we can (up to negligible errors) ignore off-diagonal entries of
$(g_{ij})$ and $(g^{ij})$ in a small neighborhood of $p$. Also,
note that $ \Big{\vert} h' (t)^2 - 1 \Big{\vert} \le \Big{\vert}
h' (t) + 1 \Big{\vert} \cdot \Big{\vert} h' (t) - 1 \Big{\vert}
\le 2 \cdot \Big{\vert} h' (t) - 1 \Big{\vert}$.
\bigskip

We will now examine the coefficients showing up in $scal$ more closely: \medskip

\noindent In the case where $a = b \ge 2$ we see that the coefficients of $h'' (t) \cdot \frac{\partial}{\partial x_1} \left(g_{\R} + g_s \right)_{aa}$
involving diagonal entries of $(g^{ij})$ are $- g^{11} \cdot g^{aa}$ satisfying $(1 + \gamma)^2 \ge \Big{\vert} - g^{11} \cdot g^{aa} \Big{\vert} \ge (1 -
\gamma)^2 >
0$.\\

On the other hand, setting $I_6$ to be the set of the remaining terms in $I_5$, those with $a \neq b$, we see that each of these $\frac{\partial}{\partial x_1}
\left(g_{\R} + g_s \right)_{ab}$ has coefficients of the form $g^{ij} \cdot g^{\alpha \beta}$ at least one of them being an off-diagonal entry:
\\ checking through the terms $g^{ij} \cdot ( \frac{\partial \Gamma^k_{ij}}{\partial x_k} - \frac{\partial
\Gamma^k_{ik}}{\partial x_j})$ in \emph{scal} we see that either $i \neq j$ or
if $i=j$ then we have from\\

\noindent {\small $ \frac{\partial \Gamma^1_{ii}}{\partial x_1} = Q (\frac{\partial g}{\partial x}) + \frac{1}{2} \sum_p g^{1 p} \cdot \left( \frac{\partial^2
g_{i p}}{\partial x_i \partial x_1} + \frac{\partial^2 g_{i p}}{\partial x_i
\partial x_1} - \frac{\partial^2 g_{i i}}{\partial x_p \partial x_1} \right) $
and $ \frac{\partial \Gamma^k_{1k}}{\partial x_1} = Q (\frac{\partial g}{\partial x}) + \frac{1}{2} \sum_p g^{k p} \cdot \left( \frac{\partial^2 g_{1
p}}{\partial x_k
\partial x_1} + \frac{\partial^2 g_{k p}}{\partial x_1 \partial x_1} - \frac{\partial^2 g_{k 1}}{\partial x_p \partial x_1} \right) $}\\

In the first case the only way would be $i=1$ but that implies $p \neq 1$, in the second one we get $k \neq p$.\\

Hence for a sufficiently small neighborhood of $p \in M$ we will have that:

\[ \Big{\vert} \sum_{I_6} \frac{\partial}{\partial x_1} \left(g_{\R} + g_s
\right)_{ab} \cdot P_{2n} \left( (g_{\R} + g_s)_{vw} \right) \cdot
\det \left( (g_{\R} + g_s)_{\alpha \beta} \right)^{-2}
\Big{\vert}_{s = h (t)} \Big{\vert} < \varepsilon \]

Also, from the positivity of the mean curvature of ${\cal{W}}^{n}$ , we get at $p$ using $(D)$: \\
$\sum_{1<m \le n} \frac{\partial g_{mm}}{\partial x_1} \vert_{t=0} < 0$ and, using a sufficiently small $\gamma$ for all $i=2,... n$ we obtain on a suitably
small neighborhood of $p$, for the diagonal terms involved in $I_5$

\[ h''(t) \cdot \sum \frac{\partial}{\partial x_1} \left(g_{\R} + g_s
\right)_{ii} \cdot P_{2n} \left( (g_{\R} + g_s)_{vw} \right) \cdot
\det \left( (g_{\R} + g_s)_{\alpha \beta} \right)^{-2}
\Big{\vert}_{s = h (t)} = \]
\[ = h''(t) \cdot \sum - \frac{\partial}{\partial x_1} \left(g_{\R} + g_s
\right)_{ii} \cdot \left(g_{\R} + g_s \right)^{11} \cdot \left(g_{\R} + g_s \right)^{ii} {\vert}_{s = h (t)} \]
\[ = -h''(t) \cdot g^{11} \cdot \sum_{i \geq 2} \frac{\partial
g_{ii}}{\partial x_1} \cdot g^{ii} \geq -h''(t) \cdot (1 \pm \gamma) \cdot \sum_{i \geq 2} \frac{\partial g_{ii}}{\partial x_1} \cdot (1 \pm \gamma) \geq
-\frac{1}{2} \cdot h''(t) \cdot tr A_{\p } > 0,\]

Now we can estimate the right hand side of $(\ast \ast \ast)$: \\
{\footnotesize \[ \ge  - 2 \cdot |h' (t) - 1| \cdot
\max\limits_{N_{h (t)}} \Big\vert \sum\limits_{I_1}
\frac{\partial}{\partial x_1} \left(g_{\R} + g_s \right)_{ab}
\cdot \frac{\partial}{\partial x_1} \left(g_{\R} + g_s
\right)_{ij} \cdot P_{3n} \left( (g_{\R} + g_s)_{vw} \right) \cdot
\det \left( (g_{\R} + g_s)_{\alpha \beta} \right)^{-3}
\Big\vert_{s = h (t)} \Big{\vert} \]
\[ - |h' (t) - 1| \cdot \max\limits_{N_{h(t)}} \Big\vert \sum_{I_2\;(k \geq 2)} \frac{\partial}{\partial x_1} \left(g_{\R} +
g_s \right)_{ab} \cdot \frac{\partial}{\partial x_k} \left(g_{\R}
+ g_s \right)_{ij} \cdot P_{3n} \left( (g_{\R} + g_s)_{vw} \right)
\cdot \det \left( (g_{\R} + g_s)_{\alpha \beta} \right)^{-3}
\Big\vert_{s = h (t)} \Big\vert \]
\[- 2 \cdot |h' (t) - 1| \cdot \max\limits_{N_{h(t)}} \Big\vert \sum\limits_{I_3} \frac{\partial^2}{\partial x_1
\partial x_1} \left(g_{\R} + g_s \right)_{ab} \cdot P_{2n} \left(
(g_{\R} + g_s)_{vw} \right) \cdot \det \left( (g_{\R} + g_s)_{\alpha \beta} \right)^{-2} \Big{\vert}_{s = h (t)} \Big\vert \]
\[ - |h' (t) - 1| \cdot \max\limits_{N_{h (t)}} \Big{\vert}
\sum\limits_{I_4\;(c \geq 2)} \frac{\partial^2}{\partial x_1
\partial x_c} \left(g_{\R} + g_s \right)_{ab} \cdot P_{2n}
\left( (g_{\R} + g_s)_{vw} \right) \cdot \det \left( (g_{\R} + g_s)_{\alpha \beta} \right)^{-2} \Big{\vert}_{s = h (t)} \Big{\vert}\]
\[ - \; h'' (t) \cdot \ve \; + \frac{1}{2} \cdot h''(t) \cdot tr A_{\p }
\ge - |h' (t) - 1| \cdot K - (\frac{1}{2} \cdot tr A_{\p } + \ve) \cdot h''(t) \mbox { for some } K>0. \]}

For sufficiently small $\ve$ we have $\ve + \frac{1}{2} \cdot tr
A_{{\cal{W}}^{n} } < 0$ and using lemma (10.3) we conclude that
the last right hand side is $>0$, proving that $scal$ remains
positive (actually gains positivity) under the deformation $g_{\R}
+ g_t \rightarrow g_{\R} + g_{h(t)}$.\qed

\end{document}